\documentclass[10pt]{amsart}
\usepackage{amssymb}
\usepackage{latexsym}
\newtheorem{prop}{Proposition}
\newtheorem{defn}{Definition}
\newtheorem{thm}{Theorem}
\newtheorem{cor}{Corollary}

\theoremstyle{remark}

\def\bP{{\mathbb P}}

\def\Z{{\mathbb Z}}

\def\C{{\mathbb C}}
\def\cO{{\mathcal O}}

\def\D{{\mathcal D}}
\def\AS{{\mathfrak S}}
\def\BS{{\mathfrak B}}
\def\CS{{\mathfrak C}}
\def\DS{{\mathfrak D}}
\def\l{{\lambda}}
\def\m{{\mu}}
\def\n{{\nu}}

\def\X{{\mathfrak X}}

\def\a{{\alpha}}
\def\b{{\beta}}
\def\g{{\gamma}}
\def\s{{\sigma}}
\def\om{{\varpi}}

\DeclareMathOperator{\Hom}{Hom}
\DeclareMathOperator{\Spec}{Spec}
\DeclareMathOperator{\rk}{rk}

\newcommand{\dis}{\displaystyle}

\newcommand{\ssm}{\smallsetminus}
\newcommand{\gequ}{\geqslant}
\newcommand{\lequ}{\leqslant}
\newcommand{\lra}{\longrightarrow}
\newcommand{\hra}{\hookrightarrow}
\newcommand{\ra}{\rightarrow}
\newcommand{\skipline}{\vspace{\baselineskip}}

\newcommand{\Ker}{\mbox{Ker}}

\newcommand{\ov}{\overline}
\newcommand{\noin}{\noindent}
\newcommand{\wt}{\widetilde}
\newcommand{\Pf}{\mbox{Pfaffian}}
\newcommand{\wh}{\widehat}

\newcommand{\vlra}{\relbar\joinrel\longrightarrow}
\newcommand{\vvlra}{\relbar\joinrel\vlra}
\newcommand{\vvvlra}{\relbar\joinrel\vvlra}
\newcommand{\vvvvlra}{\relbar\joinrel\vvvlra}
\newcommand{\vvvvvlra}{\relbar\joinrel\vvvvlra}

\begin{document}

\title[Double Schubert Polynomials and Degeneracy Loci]
{Double Schubert Polynomials and Degeneracy Loci
for the Classical Groups}
\author{Andrew Kresch and Harry Tamvakis}
\date{November 7, 2001}
\thanks{The authors were supported in part by National Science 
Foundation Postdoctoral Research Fellowships.}
\address{Department of Mathematics, University of Pennsylvania,
Philadelphia, PA 19104--6395, USA}
\email{kresch@math.upenn.edu}
\address{Department of Mathematics, Brandeis University,
Waltham, MA 02454--9110, USA}
\email{harryt@brandeis.edu}

\begin{abstract}
We propose a theory of double Schubert polynomials
$P_w(X,Y)$ for the Lie types $B$, $C$, $D$
which naturally extends the family of Lascoux
and Sch\"utzenberger in type $A$. These
polynomials satisfy
positivity, orthogonality and stability properties, and represent
the classes of Schubert varieties and degeneracy loci of vector bundles.
When $w$ is a maximal Grassmannian element of the Weyl group, 
$P_w(X,Y)$ can be expressed in terms of
Schur-type determinants and Pfaffians, in analogy with 
the type $A$ formula of Kempf and Laksov. An example,
motivated by quantum cohomology, shows there are no Chern class formulas 
for degeneracy loci of `isotropic morphisms' of bundles.
\end{abstract}

\maketitle

\setcounter{section}{-1}

\section{Introduction}

\noindent
In recent years there has been interest in finding
natural polynomials that represent the classes of Schubert
varieties and degeneracy loci of vector bundles (see \cite{Tu} and
\cite{FP} for 
expositions). 
Our aim here is to define and study
polynomials which we propose as type $B$, $C$ and
$D$ double Schubert polynomials. 
Special cases of these polynomials
provide orthogonal and symplectic analogues
of the determinantal formula of Kempf and Laksov \cite{KL}.

For the general linear group
the corresponding objects are the double Schubert
polynomials $\{\AS_{\om}(X,Y),\ \om\in S_n\}$ 
of Lascoux and Sch\"utzenberger \cite{LS} \cite{L}. These type $A$
polynomials
possess a series of remarkable properties, and it is desirable to have
a theory for the other types with as many of them as possible. 
Fomin and Kirillov \cite{FK} have shown that a theory of (single)
Schubert polynomials in types $B$, $C$ and $D$ cannot satisfy all of the
type $A$ properties simultaneously. The theory developed here has 
qualities which are desirable from both the geometric and combinatorial
points of view. Our three families of double Schubert polynomials 
have a `common core', consisting of those polynomials which correspond to 
elements of the symmetric group $S_n$, sitting inside the 
respective Weyl groups.
The polynomials in this core have positive integer coefficients
and are obtained by specializing type $A$ double Schubert polynomials.
This has a natural geometric interpretation, using the inclusion of the
corresponding isotropic flag bundle into the partial $SL_{2n}$-flag bundle
obtained by omitting the isotropicity condition.

When restricted to maximal Grassmannian elements of the Weyl group, 
the single versions of
our polynomials are the $\wt{P}$- and $\wt{Q}$-polynomials of
Pragacz and Ratajski \cite{PR}. The latter objects are polynomials 
in the Chern roots of the tautological vector bundles
over  maximal isotropic Grassmannians, which represent the Schubert 
classes in the cohomology (or Chow ring) of the base variety. 
In this sense, they play a role
in types $B$, $C$ and $D$ analogous to that of Schur's $S$-functions
in type $A$ (note that Schur's $Q$-functions, introduced in \cite{Sh} to
study projective representations of the symmetric and alternating groups, 
do not have the same geometric property \cite[6.11]{P2}). 
The utility of the $\wt{P}$- and $\wt{Q}$-polynomials in
the description of
(relative) Schubert calculus and degeneracy loci was established in \cite{PR}.
Moreover, according to
\cite{Ta} and \cite{KT}, the multiplication of 
$\wt{Q}$-polynomials describes both {\em arithmetic} and {\em quantum} 
Schubert calculus on the Lagrangian Grassmannian
(see also \cite{KT2}).
Thus, the double Schubert polynomials in this paper are closely related 
to natural families of representing polynomials.

Our motivation for this work was the search for an explicit
general formula for Lagrangian and orthogonal degeneracy loci;
to place this problem in context we 
first recall the relevant results in the setting of type $A$.
Let
$Q$ and $V$ be vector bundles of ranks 
$n$ and 
$N=m+n$ respectively 
on an algebraic variety $\X$, and consider a morphism of
vector bundles
$\psi\colon V\ra Q$.
Assume that we have a complete filtration
\[
0=V_0 \subset V_1\subset V_2 \subset \cdots \subset
V_N =V
\]
of $V$ by subbundles such that $\dim V_i=i$ for all $i$.
We are also given a rank sequence $0< r_1<r_2<\cdots < r_m \lequ N$
which corresponds to an integer partition $\l=(\l_i)$ with 
$\l_i=n+i-r_i$, and let $|\l|$ denote the sum of the parts of
$\l$. 
There is the {\em degeneracy locus}
\begin{equation}
\label{1way}
\X_{\l}'  = \{\,x\in \X\ |\ \rk(V_{r_i}(x)\stackrel{\psi}\lra 
     Q(x))
     \lequ r_i-i \ \ \ \mathrm{for} \ \mathrm{all} \  i\,\}.
\end{equation}
Assume for simplicity that $\X$ is smooth and that
$\X'_{\l}$
has codimension $|\l|$ in $\X$. 
According to \cite{KL},
the class $[\X'_{\l}]$
of this locus
in the Chow ring $CH(\X)$ is given by a Schur 
determinant:
\begin{equation}
\label{kl1}
[\X'_{\l}] =\det(c_{\l_i+j-i}(Q-V_{r_i}))_{i,j} 
\end{equation}
where each $c_k(E-F)$ is defined by
the Chern class equation $c(E-F)=c(E)c(F)^{-1}$. 
Suppose now that $\psi$ is surjective, and consider the exact
sequence
\[
0 \lra S \lra V \stackrel{\psi}\lra Q \lra 0.
\]
In this case 
$\X'_{\l}$ coincides with the locus
\begin{equation}
\label{3way}
\X_{\l} = \{\,x\in \X\ |\ \dim(S(x)\cap V_{r_i}(x))\gequ i
\ \ \ \mathrm{for}  \ \mathrm{all} \ i\,\}.
\end{equation}

We would like to consider a direct analogue of (\ref{1way}) in
types $B$, $C$ and $D$ (see Section \ref{lagdl} for the 
precise definition). Unfortunately, as we shall see by
example (Section \ref{example}), there are no Chern class
formulas for degeneracy loci in this level of generality.
We consider instead the analogue of the locus (\ref{3way})
for the other types, and begin with the Lagrangian case.
Here $V$ is a rank $2n$ symplectic vector bundle over $\X$,
and $E$ and $F$ are
Lagrangian subbundles of $V$. The bundle $F$ comes with a complete
filtration $F_{\bullet}$ by subbundles $F_i\subset F$, $1\lequ i\lequ n$.
For each (strict) partition $\l$ with $\ell$ nonzero parts
the degeneracy locus
$\X_{\l}\subset \X$ is defined by
\[
\X_{\l} = \{\,x\in \X\ |\ \dim(E(x)\cap F_{n+1-\l_i}(x))\gequ i
\ \ \ \mathrm{for}  \ \ 1\lequ i\lequ \ell \,\}.
\]
Assuming that $\X_{\l}$ has
codimension $|\l|$ in $\X$, our main geometric 
result (Corollary \ref{lagcor}) is a determinantal formula for the
the class $[\X_{\l}]\in CH^{|\l|}(\X)$ 
as a polynomial in the Chern classes of
the bundles $E$ and $F_i$, which is a type $C$ analogue of (\ref{kl1}).
Corollaries \ref{lagcorB} and \ref{lagcorD} 
solve the analogous problem in the two orthogonal cases. 
These results answer a question of Fulton and Pragacz \cite[Section 9.5]{FP}; 
note that
a priori it is not clear why $[\X_{\l}]$ should be expressed as {\em
any} polynomial in the Chern classes of the bundles. Our formulas
generalize those obtained by Pragacz and Ratajski \cite{PR} for some
special cases of these loci. 

The main ingredients used in the proofs are the geometric 
work of Fulton \cite{F2} \cite{F3} and Graham \cite{Gra} and 
the algebraic tools developed by Lascoux, Pragacz and
Ratajski \cite{PR} \cite{LP1} \cite{LP2}. 
The degeneracy locus formula (\ref{kl1}) for
$\X_{\l}'$ was generalized to maps of
flagged vector bundles by Fulton \cite{F1}, using
the work of Bernstein, Gelfand, Gelfand \cite{BGG} and 
Demazure \cite{D1} \cite{D2} and 
the double Schubert polynomials of 
Lascoux and Sch\"utzenberger \cite{LS} \cite{L}. Fulton later
extended the geometric part of this story to the other classical
groups \cite{F2} \cite{F3}; however there 
is no general degeneracy
locus formula for {\em morphisms} between bundles 
in types $B$, $C$ and $D$ (note that there are such formulas
for morphisms with symmetries; see 
\cite{HT} \cite{JLP} \cite{P} \cite{F3} \cite{PR} \cite{LP3}).

To elaborate further, let $V$ be a vector bundle on $\X$ equipped with
a nondegenerate symplectic or orthogonal form.
For each element $w$ in the corresponding 
Weyl group there is a degeneracy locus $\X_w$, defined using
the attitude of an isotropic flag $E_{\bullet}$ with respect
to a fixed complete isotropic flag $F_{\bullet}$ of subbundles 
of $V$ (see \cite{F3} \cite{FP} or Sections 2, 3). 
In general one has an algorithm to write $[\X_w]$ as
a polynomial $P_w$ in the first Chern classes $X$, $Y$
of the quotient line
bundles of the flags, by applying divided difference operators 
to a kernel $P_{w_0}(X,Y)$ (which corresponds to the longest Weyl group
element $w_0$). The polynomial $P_{w_0}$ represents the class of the 
diagonal in the flag bundle of $V$, and 
is defined only modulo an ideal of relations.
Each specific choice of $P_{w_0}$ produces a different
family of representing polynomials $P_w$, which are candidates
for {\em double Schubert polynomials}. The above construction ensures that
the resulting theory of polynomials has direct geometric significance.

 As explained earlier, unlike the situation in type
$A$, it is no longer clear which kernel will give the most
desirable theory. For (single)
Schubert polynomials there have been several works in this
direction \cite{BH} \cite{FK} \cite{PR} \cite{LP1} \cite{LP2}.
Candidates for double Schubert polynomials in types $B$, $C$, $D$ 
are implicit in the works \cite{F2} \cite{F3} and 
\cite[Appendix B3]{LP1}. The double Schubert polynomials proposed
in this article differ from both of these sources.
We choose a kernel $P_{w_0}(X,Y)$ leading to
polynomials $P_w(X,Y)$ (where $P\in\{\BS, \CS, \DS\}$ depends on
the Lie type) which have the following two main properties:

\medskip
\noin
(i) 
(Positivity) When $w$ is a permutation in the symmetric group $S_n$,
$P_w(X,Y)$ is equal (in types $B$ and $C$) or closely 
related (in type $D$)
to the type $A$ double Schubert 
polynomial $\AS_{\om_0 w \om_0}(\om_0X,-Y)$,
where $\om_0$ denotes the permutation of longest length in $S_n$. 
In particular, $\BS_w(X,Y)$ and $\CS_w(X,Y)$ have 
positive integer coefficients.

\medskip
\noin
(ii) (Maximal Grassmannian)
When $w$ is a maximal Grassmannian element of the Weyl group, 
$P_w(X,Y)$ can 
be expressed by an explicit formula involving Schur-type
determinants and Pfaffians.

\medskip

In type $C$ our polynomials
specialize (when $Y=0$)
to the {\em symplectic Schubert polynomials} $\CS_w(X)$ 
of \cite{PR} \cite{LP1}, and, for maximal Grassmannian elements $w$, 
to the $\wt{Q}$-polynomials of \cite{PR}. Properties (i), (ii) above
combined with the results of \cite{PR} \cite{LP1} \cite{LP2}
can be used to push the theory further, and obtain

\medskip
\noin
(iii) (Orthogonality) Products of the single Schubert polynomials
$P_w(X)=P_w(X,0)$ give a (positive coefficient)
orthonormal basis for the full polynomial ring, as a 
module over the ring of Weyl group invariants. This product basis
corresponds in geometry to a split form of the fibration of the isotropic
flag bundle over the maximal Grassmannian bundle, 
with fibers isomorphic to the type $A$ flag variety (see 
Table \ref{table} in Section \ref{sc}).

\medskip
\noin
(iv) (Stability)
The single Schubert 
polynomials $P_w(X)$ satisfy a stability property under the
natural inclusions of the Weyl groups. The double Schubert polynomials
$P_w(X,Y)$ are stable for certain special Weyl group elements,
related to the loci for morphisms with symmetries referred to earlier.

\medskip

The hyperoctahedral group is a semidirect product of $S_n$ with 
$(\Z_2)^n$. A central theme of this work is that the 
restrictions of double Schubert polynomials in types $B$, $C$, $D$ 
to the two factors in this product (with the second realized by the
maximal Grassmannian elements) are amenable to study and should be
related to previously known families of polynomials. In fact, our
polynomials are
characterized (Proposition \ref{char}) by specifying the values
of their single versions for 
the maximal elements in each factor, assuming the kernel $P_{w_0}$ 
satisfies a `Cauchy formula' (Corollary \ref{Cauchy}).

Our interest in these questions originated in an application to quantum
cohomology. In the $SL_N$ case, Bertram
\cite{bertram} used (\ref{kl1}) to prove a `quantum Giambelli
formula' for the Schubert classes in the small quantum cohomology 
ring of the Grassmannian. The degeneracy locus problem in loc.\
cit.\ occurs on a Quot scheme \cite{grothendieck}
and it is crucial that (\ref{kl1}) holds for an arbitrary morphism
$\psi$. There is a type $C$ analogue of the Quot scheme,
and one can ask whether natural Lagrangian analogues of 
$[\X'_{\l}]$ can be expressed in terms of the 
Chern classes of the bundles involved. In Section \ref{example}
we show that
{\em no such formula exists} (Proposition \ref{nopoly})
by analyzing the structure of 
the Quot scheme
$LQ_1(2,4)$, which compactifies the space of degree $1$ maps
from $\bP^1$ to the Lagrangian Grassmannian $LG(2,4)$. We hope
this example is of independent interest.

Here is a brief outline of this article. In Section \ref{A} we
introduce the type $C$ double Schubert polynomials $\CS_w(X,Y)$
and prove their basic algebraic properties.
Section \ref{dl} connects this work to the
geometry of symplectic and Lagrangian
degeneracy loci. The double Schubert polynomials and
loci for the orthogonal 
types $B$ and $D$ are studied in Section \ref{C}. Finally, 
Section \ref{example} presents the example of the Lagrangian
Quot scheme $LQ_1(2,4)$. 

The authors would like to thank
Piotr Pragacz for communication regarding his 
work with Lascoux \cite{LP1} \cite{LP2}. The computational techniques 
with symplectic and orthogonal divided difference operators
developed in these papers are essential tools in the present work.
We also wish to thank
Sara Billey for her comments on an earlier version of
this paper, and an anonymous referee for remarks which prompted us to
develop our theory of polynomials systematically. 
It is a pleasure to thank the Institut des Hautes \'Etudes 
Scientifiques for its stimulating atmosphere and hospitality in the
fall of 2000.

\section{Type $C$ double Schubert polynomials}
\label{A}

\subsection{Initial definitions}
\label{definitions}
Let us begin with some combinatorial preliminaries: for the most
part we will follow the notational conventions of \cite[Section 1]{M3};
although not strictly necessary, the reader may find it helpful to
identify an integer partition $\l=(\l_1,\ldots,\l_r)$ with 
its Young diagram of boxes. The sum $\sum \l_i$ of the
parts of $\l$ is the {\em weight} $|\l|$ and the number of (nonzero)
parts is the {\em length} $\ell(\l)$ of $\l$. We set $\l_r=0$ for any
$r>\ell(\l)$.
Define the containment relation 
$\m\subset \l$ for partitions by the inclusion of their
respective diagrams. The union $\l\cup \m$, intersection $\l\cap\m$
and set-theoretic difference $\l\ssm\m$ are
defined using the (multi)sets of parts of $\l$ and $\m$. 
A partition is {\em strict}
if all its nonzero parts are different; we define $\rho_n=
(n,n-1,\ldots,1)$ and let $\D_n$ be the set of strict partitions
$\l$ with $\l\subset \rho_n$.
For $\l\in \D_n$, the {\em dual} partition $\l'=\rho_n\ssm\l$ 
is the strict partition
whose parts complement the parts of $\l$ in the set $\{1,\ldots,n\}$.

Define the {\em excess} $e(\l)$ of a strict partition
$\l$ by
\[
e(\l)=|\l|-(1+\cdots + \ell(\l))=|\l|-\ell(\l)(\ell(\l)+1)/2.
\]
More generally, given $\a, \b \in \D_n$ with $\a\cap\b=\emptyset$,
define an {\em intertwining number} 
$e(\a,\b)$ as follows: for each $i$ set
\[
m_i(\a,\b)=\#\,\{\,j\ |\ \a_i>\b_j>\a_{i+1} \,\}
\]
and
\[
e(\a,\b)=\sum_{i\gequ 1} i\,m_i(\a,\b).
\]
Note that for any $\l\in \D_n$ we have $e(\l,\l')=e(\l)$.

We will use multiindex notation for sequences of commuting independent
variables;
in particular for any $k$ with $1\lequ k\lequ n$ let
$X_k=(x_1,\ldots,x_k)$ and $Y_k=(y_1,\ldots,y_k)$; also set
$X=X_n$ and $Y=Y_n$. Following Pragacz and Ratajski \cite{PR}, 
for each
partition $\l$ with $\l_1\lequ n$, we define a symmetric polynomial 
$\wt{Q}_{\l}\in \Z[X]$ as follows: set $\wt{Q}_i(X)=e_i(X)$ to be the
$i$-th elementary symmetric polynomial in the variables $X$.
For $i,j$ nonnegative integers let
\[
\wt{Q}_{i,j}(X)=
\wt{Q}_i(X)\wt{Q}_j(X)+
2\sum_{k=1}^j(-1)^k\wt{Q}_{i+k}(X)\wt{Q}_{j-k}(X).
\]
If $\l=(\l_1\gequ\l_2\gequ\cdots\gequ\l_r\gequ 0)$ is a partition with
$r$ even (by putting $\l_r=0$ if necessary), set
\[
\dis
\wt{Q}_{\l}(X)=\Pf[\wt{Q}_{\l_i,\l_j}(X)]_{1\lequ i<j\lequ r}.
\]
These $\wt{Q}$-polynomials 
are modelled on Schur's $Q$-polynomials \cite{Sh}
and were used in \cite{PR} to describe certain 
Lagrangian and orthogonal degeneracy loci; we will generalize these
results in the following sections. We also need
the reproducing kernel
\[
\wt{Q}(X,Y)=\sum_{\l\in \D_n}\wt{Q}_{\l}(X)\wt{Q}_{\l'}(Y).
\]
For more information on these polynomials we refer to 
\cite{PR} \cite{LP1}.

Let $S_n$ denote the symmetric group of permutations of the set
$\{1,\ldots,n\}$ whose elements $\s$ are written in single-line
notation as $(\s(1),\s(2),\ldots,\s(n))$ (as usual we will write all
mappings on the left of their arguments.) $S_n$ is the Weyl group
for the root system $A_{n-1}$ and is generated by the simple
transpositions $s_i$ for $1\lequ i\lequ n-1$, where $s_i$ 
interchanges $i$ and $i+1$ and fixes all other elements of $\{1,\ldots,n\}$.

The hyperoctahedral group $W_n$ is the Weyl group for the root
systems $B_n$ and $C_n$; the elements of $W_n$ are permutations
with a sign attached to each entry. We will adopt the notation
where a bar is written over an element with a negative sign 
(as in \cite[Section 1]{LP1}); the latter are the
{\em barred} entries, and their positive counterparts are 
{\em unbarred}.
$W_n$ is an extension of $S_n$ by an
element $s_0$ which acts on the right by
\[
(u_1,u_2,\ldots,u_n)s_0=(\ov{u}_1,u_2,\ldots,u_n).
\]
The elements of maximal length in $S_n$ and $W_n$ are
\[
\om_0=(n,n-1,\ldots,1) \  \ \ \mathrm{and} \ \ \ 
w_0=(\ov{1},\ov{2},\ldots,\ov{n})
\]
respectively. Let $\l\subset \rho_n$ be a strict partition,
$\ell=\ell(\l)$ and $k=n-\ell=\ell(\l')$. The barred
permutation 
\[
w_{\l}=(\ov{\l}_1,\ldots,\ov{\l}_{\ell},\l'_k,\ldots,\l'_1)
\]
is the {\em maximal Grassmannian element} of $W_n$ corresponding
to $\l$ (for details, see \cite[Prop.\ 1.7]{LP1} and 
\cite{BGG} \cite{D1} \cite{D2}).

The group $W_n$ acts on the ring $A[X]$ of polynomials in
$X$ with coefficients in any commutative ring $A$:
the transposition $s_i$ interchanges $x_i$ and $x_{i+1}$
for $1\lequ i\lequ n-1$, while $s_0$ replaces $x_1$ by $-x_1$
(all other variables remain fixed). The ring of invariants 
$A[X]^{W_n}$ is the ring of polynomials in $A[X]$ symmetric
in the $X^2=(x_1^2,\ldots,x_n^2)$.

Assume that $2$ is not a zero divisor in $A$.
Following \cite{BGG} and 
\cite{D1} \cite{D2}, there are divided difference
operators $\partial_i\colon  A[X]\ra A[X]$. For $1\lequ i\lequ n-1$
they are defined by
\[
\partial_i(f)=(f-s_if)/(x_i-x_{i+1})
\]
while 
\[
\partial_0(f)=(f-s_0f)/(2x_1),
\]
for any $f\in A[X]$. For each $w\in W_n$, define an operator
$\partial_w$ by setting
\[
\partial_w=\partial_{a_1}\circ \cdots \circ \partial_{a_r}
\]
if $w=s_{a_1}\cdots s_{a_r}$ is a reduced decomposition of 
$w$ (so $r=\ell(w)$). One can show that $\partial_w$ is
well-defined, using the fact that the operators $\partial_i$
satisfy the same set of Coxeter relations as the generators
$s_i$, $0\lequ i\lequ n-1$. In the applications here and in
the next section we will
take $A=\Z[Y]$. For each $i$ with $1\lequ i\lequ n-1$ we
set $\partial'_i=-\partial_i$, while $\partial'_0:=\partial_0$.
For each $w\in W_n$ we use $\partial'_w$ to denote the 
divided difference operator where we use the $\partial'_i$'s
instead of the $\partial_i$'s. 

Following Lascoux and Sch\"utzenberger \cite{LS} \cite{L} 
we define, for each $\om\in S_n$,
a type $A$
{\em double Schubert polynomial} $\AS_{\om}(X,Y)\in\Z[X,Y]$ by
\[
\AS_{\om}(X,Y)=\partial_{\om^{-1}\om_0}\left(\Delta(X,Y)\right)
\]
where 
\[
\Delta(X,Y)=\prod_{i+j\lequ n}(x_i-y_j).
\]
Our main references for these polynomials are \cite{M1}
\cite{M2}. The type A (single) {\em Schubert polynomial} is 
$\AS_{\om}(X):=\AS_{\om}(X,0)$. The double and single polynomials
are related by the formula
\cite[(6.3)]{M2}
\begin{equation}
\label{dtos}
\AS_{\om}(X,Y)=\sum_{u,v}\AS_u(X)\AS_v(-Y),
\end{equation}
summed over all $u,v\in S_n$ such that $u=v\om$ and 
$\ell(\om)=\ell(u)+\ell(v)$.

\subsection{Main theorems}
\label{mts}
We intend to give analogues of the polynomials $\AS_{\om}(X,Y)$
for the other classical groups, and begin with the symplectic group.
Our definition is motivated by the properties that follow and
the applications to the geometry of degeneracy loci in 
Section \ref{dl}.

\begin{defn}
\label{Cdefn}
For every $w\in W_n$ the {\em type $C$ 
double Schubert polynomial $\CS_w(X,Y)$} is given by
\[
\CS_w(X,Y)= (-1)^{n(n-1)/2}\partial'_{w^{-1}w_0}\left(
\Delta(X,Y)\wt{Q}(X,Y)\right).
\]
\end{defn}

The polynomials $\CS_w(X):=\CS_w(X,0)$ have already been
defined by Lascoux, Pragacz and Ratajski \cite{PR} \cite[Appendix
A]{LP1}, where they are called {\em symplectic Schubert polynomials}.
It is shown in loc.\ cit.\ that the $\CS_w(X)$ enjoy `maximal 
Grassmannian' and `stability' properties. We will see that the type
$C$ double Schubert polynomials $\CS_w(X,Y)$
satisfy an analogue of the former
but are not stable under the natural embedding $W_{n-1} \hra W_n$
in general. Note that the $\CS_w(X,Y)$
differ from the polynomials found in 
\cite[Appendix B3]{LP1}, because the divided difference operators in
loc.\ cit.\ (which agree with those in \cite{F3}) 
differ from the ones used in Definition \ref{Cdefn}. 
\footnote{Lascoux and Pragacz have informed us that the polynomials 
$\CS_w(X,Y)$ were 
defined in a preliminary version of [LP1].}
For the precise connection of the $\CS_w(X,Y)$ with geometry,
see the proof of Theorem \ref{lagfulthm}.

Let $\CS_{\l}(X,Y)=\CS_{w_{\l}}(X,Y)$ for each partition $\l\in \D_n$.
In type $A$, when $\om\in S_n$ is a Grassmannian permutation $\om_{\l}$,
the Schubert polynomial $\AS_{\om_{\l}}(X,Y)$ is a {\em multi-Schur
function} (see \cite[(6.14)]{M2}). We may thus regard the $\CS_{\l}(X,Y)$
as type $C$ analogues of multi-Schur determinants. Our first theorem
gives a determinantal expression for these functions.

\begin{thm}[Maximal Grassmannian]
\label{mgras}
For any strict partition
$\l\in \D_n$ the double Schubert polynomial
$\CS_{\l}(X,Y)$ is equal to 
\[
(-1)^{e(\l)+|\l'|}
\sum_{\a}\wt{Q}_{\a}(X)
\sum_{\b}
(-1)^{e(\a,\b)+|\b|}
\wt{Q}_{(\a\cup \b)'}(Y)
\det(e_{\b_i-\l'_j}(Y_{n-\l'_j})),
\]
where the first sum is over all $\a\in \D_n$ and the
second over $\b\in \D_n$ with $\b\supset \l'$,
$\ell({\b})=\ell({\l'})$ and $\a\cap\b=\emptyset$.
\end{thm}

\begin{proof} We argue along the lines of the proof of
\cite[Theorem A.6]{LP1}.
Let $\ell=\ell(\l)$, $k=n-\ell$ and 
observe that $w_{\l}=w_0\,\tau_{\l}\,\om_k\,\delta_k^{-1}$, where 
\begin{gather*}
\tau_{\l} = (\l'_k,\ldots,\l'_1,\l_1,\ldots,\l_{\ell}) \\
\om_k = (k,\ldots ,2,1,k+1,\ldots ,n) \\
\delta_k = (s_{\ell}\cdots s_1s_0)(s_{\ell +1}\cdots s_1s_0)
\cdots (s_{n-1}\cdots s_1s_0);
\end{gather*}
it follows that
\[
\partial_{w_{\l}^{-1}w_0}=\partial_{\delta_k}\circ
\partial_{\om_k}\circ\partial_{\tau_{\l}^{-1}}.
\]
Now compute
\begin{align}
\label{lcalc}
\partial_{\tau_{\l}^{-1}}(\Delta(X,Y)) &=
\AS_{\om_0\tau_{\l}}(X,Y) \\
\label{e3}
&= \prod_{j=1}^k\prod_{p=1}^{n-\l'_j}(x_{k+1-j}-y_p) \\
\label{e4}
&= \sum_{\gamma}(-1)^{|\gamma|-|\l'|}
\prod_{j=1}^k x_{k+1-j}^{n-\gamma_j}e_{\gamma_j-\l'_j}
(Y_{n-\l'_j}).
\end{align}
In the above calculation equality (\ref{e3})
holds because the permutation 
\[
\om_0\tau_{\l}=(n+1-\l'_k,\ldots,n+1-\l'_1,n+1-\l_1,\ldots,
n+1-\l_{\ell})
\]
is {\em dominant} \cite[(6.14)]{M2}, and 
the sum (\ref{e4}) is over all $k$-tuples $\gamma=(\gamma_1,
\ldots,\gamma_k)$ of nonnegative integers.

Observe that 
\[
\partial'_{\delta_k}\circ\partial'_{\om_k}=(\nabla_1)^k
\]
where $\nabla_1=\partial'_{n-1}\cdots\partial_1'\partial_0$.
The action
of $(\nabla_1)^k$ relevant to our computation is explained by
Lascoux and Pragacz; if $\n_i\lequ n-i$
for $1\lequ i\lequ k$ then \cite[Lemma 5.10]{LP1}:
\[
(\nabla_1)^k(x_1^{\n_1}\cdots x_k^{\n_k}\cdot f)=
\nabla_1(x_1^{\n_k}\cdots \nabla_1(x_1^{\n_2}\cdot
\nabla_1(x_1^{\n_1}\cdot f))\cdots )
\]
for any function $f$ symmetric in the $X$ variables. Moreover
\cite[Theorem 5.1]{LP1} states that if $\l\in \D_n$ and
$0<r\lequ n$ then $\nabla_1(x_1^{n-r}\wt{Q}_{\l}(X))$ vanishes
unless $\l_p=r$ for some $p$, in which case
\[
\nabla_1(x_1^{n-r}\wt{Q}_{\l}(X))=(-1)^{p-1}\wt{Q}_{\l\ssm r}(X).
\]

Using the above facts and equation (\ref{e4}) we deduce that
\begin{gather*}
\partial'_{w_{\l}^{-1}w_0}(\Delta(X,Y)\, \wt{Q}(X,Y))=
(\nabla_1)^k\left(\wt{Q}(X,Y) 
\partial'_{\tau_{\l}^{-1}}(\Delta(X,Y))\right) \\
= (\nabla_1)^k\Bigg(\wt{Q}(X,Y) 
\sum_{\substack{\gamma\in \D_n \\ \ell(\gamma)=k}}
(-1)^{|\gamma|-|\l'|+\ell(\tau_{\l})}
\sum_{\s\in S_k}\prod_{j=1}^k x_{k+1-j}^{n-\gamma_{\s(j)}}
e_{\gamma_{\s(j)}-\l'_j} (Y_{n-\l'_j})\Bigg) \\
= (-1)^{r(\l)}\sum_{\a}\wt{Q}_{\a}(X)
\sum_{\b}(-1)^{e(\a,\b)+|\b|}
\wt{Q}_{(\a\cup \b)'}(Y)
\det(e_{\b_i-\l'_j}(Y_{n-\l'_j}))
\end{gather*}
where the ranges of summation are as in the statement of the
theorem and $r(\l)=\ell(\tau_{\l})+k(k-1)/2+|\l'|$. 
Finally, note that
\[
\ell(\tau_{\l})=\binom{n+1}{2}+\binom{\ell}{2}-
\binom{k+1}{2} -|\l|
\]
and the signs fit to complete the proof. 
\end{proof}

\medskip
\noindent
In Section \ref{dl} we apply Theorem \ref{mgras} to obtain formulas
for Lagrangian degeneracy loci. 

\medskip

Let $*\colon S_n\ra S_n$ be the length preserving involution defined by
\[
\om^*=\om_0\om\om_0.
\]
The next theorem shows that for $\om\in S_n$, the polynomial
$\CS_{\om}(X,Y)$ is closely related to 
a type $A$ double Schubert polynomial,
with the $X$ variables in reverse order. See Section \ref{inter}
for a geometric interpretation.

\begin{thm}[Positivity]
\label{posit}
For every $\om\in S_n$,
\[
\CS_{\om}(X,Y)=\AS_{\om^*}(\om_0X,-Y).
\]
In particular, $\CS_{\om}(X,Y)$ has nonnegative integer coefficients.
\end{thm}
\begin{proof}
Let $v_0=\om_0 w_0$; then the analysis in 
\cite[Sect.\ 4]{LP1} gives
\begin{equation}
\label{dvaction}
\partial'_{v_0}(\AS_u(X)\wt{Q}(X,Y))=
(-1)^{\ell(u)}\AS_u(\om_0 X)
\end{equation}
for every $u\in S_n$. We now use (\ref{dtos}) and (\ref{dvaction})
to compute
\begin{align*}
\CS_{\om_0}(X,Y) &= (-1)^{n(n-1)/2}\partial'_{v_0}
(\Delta(X,Y)\wt{Q}(X,Y)) \\
&= (-1)^{\ell(\om_0)}\partial'_{v_0}
\left(\sum_{u,v}\AS_u(X)\AS_v(-Y)\wt{Q}(X,Y)\right) \\
&= \sum_{u,v}(-1)^{\ell(v)}\AS_u(\om_0X)\AS_v(-Y) \\
&= \AS_{\om_0}(\om_0X,-Y),
\end{align*}
where the above sums are over all $u,v\in S_n$ with $u=v\om_0$
and $\ell(u)+\ell(v)=\ell(\om_0)$.
It follows that for any $\om\in S_n$,
\begin{align}
\label{interc}
\CS_{\om}(X,Y) &= \partial'_{\om^{-1}\om_0}(\om_0\AS_{\om_0}(X,-Y)) \\
\label{intereq}
&= \om_0\partial_{\om_0\om^{-1}}\sum_{u\in S_n}\AS_u(X) \AS_{u\om_0}(Y).
\end{align}
Applying  \cite[(4.3)]{M2} gives
\begin{equation}
\label{brac}
\partial_{\om_0\om^{-1}}\AS_u(X)=\left\{ \begin{array}{cl}
       \AS_{u\om\om_0}(X) & \mathrm{ if } \,\ \ell(u)+\ell(\om)+\ell(u\om)
                                               =2\ell(\om_0), \\
       0 & \mathrm{ otherwise. }
       \end{array} \right.
\end{equation}
We now use (\ref{brac}) in (\ref{intereq}) and apply (\ref{dtos}) 
to obtain
\[
\CS_{\om}(X,Y) = \om_0\AS_{\om_0\om\om_0}(X,-Y),
\]
which is the desired result.
\end{proof}

\medskip

Combining Theorem \ref{posit} with the known results for type $A$ double
Schubert polynomials gives corresponding ones for the $\CS_{\om}(X,Y)$. For
example, we get

\begin{cor}
\label{CtoS}
For every $\om\in S_n$ we have $\CS_{\om}(X)=\AS_{\om^*}(\om_0X)$.
Moreover
\[
\CS_{\om}(X,Y)=\sum_{u,v}\CS_u(X)\CS_v(\om_0 Y)
\]
summed over all $u,v\in S_n$ with $u=v\om^*$ and $\ell(u)+\ell(v)=
\ell(\om)$.
\end{cor}

\subsection{Product basis and orthogonality}
Following Lascoux and Pragacz \cite[Sect.\ 1]{LP1} (up to a sign!),
we define a 
$\Z[X]^{W_n}$-linear scalar product
\begin{equation}
\label{inner}
\left<\ ,\ \right> : \Z[X]\times \Z[X] \lra \Z[X]^{W_n}
\end{equation}
by
\begin{equation}
\label{innerprod}
\left<f,g\right> = \partial_{w_0}(fg)
\end{equation}
for any $f,g\in \Z[X]$. 
For each partition $\l\in\D_n$ let 
$\CS_{\l}(X)=\CS_{\l}(X,0)=\wt{Q}_{\l}(X)$. We deduce from
Corollary \ref{CtoS} and \cite[(4.11)]{M1} that the set
$\{\CS_{\om}(X)\}_{\om\in S_n}$ is a free $\Z$-basis for the
additive subgroup of $\Z[X]$ spanned by the monomials
$\{x_1^{\a_1}\cdots x_n^{\a_n}\}$, with $\a_i\lequ i-1$ for each $i$.
It follows that the $\{\CS_{\om}(X)\}_{\om\in S_n}$ form a
basis of $\Z[X]$ as a $\Z[X]^{S_n}$-module.
In addition, we know from \cite[Prop.\ 4.7]{PR} that 
$\{\CS_{\l}(X)\}_{\l\in\D_n}$ is a basis of $\Z[X]^{S_n}$ as a
$\Z[X]^{W_n}$-module.

\begin{prop}[Orthogonality]
\label{Corth} {\em (i)}
The products $\{\CS_{\om}(X)\CS_{\l}(X)\}$ for $\om\in S_n$ and $\l\in\D_n$
form a basis for the
polynomial ring $\Z[X]$ as a $\Z[X]^{W_n}$-module. Moreover, we have
the orthogonality relation 
\begin{equation}
\label{orth1}
\Bigl<\CS_u(X)\CS_{\l}(X),\CS_{v\om_0}(-\om_0X)\CS_{\mu'}(X)\Bigr>=
\delta_{u,v}\delta_{\l,\mu}
\end{equation}
for all $u,v \in S_n$ and $\l,\mu\in\D_n$.

\medskip
\noin
{\em (ii)} Let $\CS_{\om,\l}(X)=\CS_{\om}(X)\CS_{\l}(X)$ and suppose 
$u,v\in S_n$ and $\l,\mu\in\D_n$ are such that $\ell(u)+\ell(v)=n(n-1)/2$.
Then
\begin{equation}
\label{orth2}
\Bigl<\CS_{u,\l}(X),\CS_{v,\mu}(X)\Bigr>=\left\{ \begin{array}{cl}
             1 &  \mathrm{ if } \ \ v=\om_0u \ \ \mathrm{and} \ \ 
                    \mu=\l', \\
             0 &  \mathrm{ otherwise. }
       \end{array} \right.
\end{equation}
\end{prop}
\begin{proof}
It is shown in \cite[Sect.\ 1]{LP1} that $\Z[X]$ is a 
free $\Z[X]^{W_n}$-module with basis $\{\AS_{\om}(X)\wt{Q}_{\l}(X)\}$,
for $\om\in S_n$, $\l\in\D_n$. Furthermore, this basis has the
orthogonality property \cite[p.\ 12]{LP1}
\[
\Bigl<\AS_u(X)\wt{Q}_{\l}(X),\AS_{v\om_0}(-\om_0X)\wt{Q}_{\mu'}(X)\Bigr>=
(-1)^{\ell(\om_0)}\delta_{u,v}\delta_{\l,\mu}.
\]
Part (i) follows immediately by applying Corollary \ref{CtoS}. To
prove (ii), factor 
$\partial_{w_0}=\partial_{v_0}\partial_{\om_0}$ (where $v_0=w_0\om_0$) 
and observe that
\begin{equation}
\label{innerfact}
\Bigl<\CS_{u,\l}(X),\CS_{v,\mu}(X)\Bigr>=
\partial_{\om_0}(\CS_u(X)\CS_v(X))\cdot 
\partial_{v_0}(\CS_{\l}(X)\CS_{\mu}(X)).
\end{equation}
Now use (\ref{innerfact})  and the 
corresponding properties of the inner products defined by $\partial_{\om_0}$
and $\partial_{v_0}$ (stated in \cite[(5.4)]{M1} and \cite[(10)]{LP1}, 
respectively).
\end{proof}

\medskip
Note that the polynomials in the basis $\{\CS_{\om,\l}(X)\}$ 
have positive coefficients, but do not 
represent the Schubert cycles in the complete symplectic flag variety
$Sp(2n)/B$. The orthogonality relations (\ref{orth1}) and (\ref{orth2})
are however direct
analogues of the ones for type $A$ Schubert polynomials
\cite[(5.4) (5.5) (5.8)]{M1}. The divided difference operator in 
(\ref{innerprod}) for the maximal length element
corresponds to a Gysin homomorphism in geometry (see \cite{BGG}
\cite{D1} \cite{AC}).

Let
$\{\CS^{\om,\l}(X)\}_{(\om,\l)\in S_n\times \D_n}$ be the 
$\Z[X]^{W_n}$-basis of $\Z[X]$ adjoint to the basis $\{\CS_{\om,\l}(X)\}$
relative to the scalar product (\ref{inner}). By Proposition \ref{Corth}
we have
\begin{equation}
\label{adjoint}
\CS^{\om,\l}(X)=\CS_{\om\om_0}(-\om_0 X)\CS_{\l'}(X).
\end{equation}
We now obtain type $C$ analogues of \cite[(5.7) and (5.9)]{M2}:

\begin{cor}[Cauchy formula]
\label{Cauchy}
We have
\begin{align*}
\CS_{w_0}(X,Y) &= \sum_{\om,\l}
\CS_{\om\om_0}(-\om_0 X)\CS_{\l'}(X)\CS_{\om}(\om_0 Y)\CS_{\l}(Y) \\
&= \sum_{\om,\l}
\CS^{\om,\l}(X)\CS_{\om,\l}(\om_0 Y),
\end{align*}
summed over all $\om\in S_n$ and $\l\in \D_n$.
\end{cor}
\begin{proof}
We use (\ref{dtos}) and the definitions in Section \ref{definitions}
to get
\begin{align*}
\CS_{w_0}(X,Y) &= \AS_{\om_0}(Y,X)\wt{Q}(X,Y) \\
&= \sum_{\om\in S_n} \AS_{\om}(Y) \AS_{\om\om_0}(-X)
\sum_{\l\in \D_n} \wt{Q}_{\l}(X)\wt{Q}_{\l'}(Y).
\end{align*}
The result now follows from  Corollary \ref{CtoS} and (\ref{adjoint}).
\end{proof}

\subsection{Examples}
1) The type $C$ double Schubert polynomials for $n=2$ are
\begin{gather*}
\CS_{(\ov{1},\ov{2})}=(y_1-x_1)\wt{Q}(X_2,Y_2)=
(y_1-x_1)(x_1+x_2+y_1+y_2)(x_1x_2+y_1y_2) \\
\CS_{(1,\ov{2})}=-x_1^2x_2+(x_2^2-y_2^2)y_1+x_2y_1^2
\ \ \ \ \ \ \ \ \ \ \ \ \ \ \ \ \ \ 
\CS_{(\ov{2},\ov{1})}=\wt{Q}(X_2,Y_2) \\
\CS_{(\ov{2},1)}=x_1x_2+(x_1+x_2)y_1+y_1^2
\ \ \ \ \ \ \ \ 
\CS_{(2,\ov{1})}=x_2^2+x_2(y_1+y_2)+y_1y_2 \\
\CS_{(2,1)}=x_2+y_1 \ \ \ \ \ \ \ \ \ \ \ \ \ \ \ \ \ 
\CS_{(\ov{1},2)}=x_1+x_2+y_1+y_2 \\
\CS_{(1,2)}=1.
\end{gather*}

\medskip
\noin
2) The type $C$ double Schubert polynomials $\CS_{\om}$ for $n=3$ 
and $\om \in S_3$ are
\begin{gather*}
\CS_{(3,2,1)}=(x_2+y_1)(x_3+y_1)(x_3+y_2) \\
\CS_{(2,3,1)}=(x_3+y_1)(x_3+y_2)
\ \ \ \ \ \ \ \ 
\CS_{(3,1,2)}=(x_2+y_1)(x_3+y_1) \\
\hspace{-0.7cm}
\CS_{(2,1,3)}=x_2+x_3+y_1+y_2 \qquad
\ \ \ \ \ \ \ \ \ \ \ 
\CS_{(1,3,2)}=x_3+y_1 \\
\CS_{(1,2,3)}=1.
\end{gather*}

\medskip
\noin
3) For fixed $n$ and $i=0,1,\ldots,n-1$ we have
\[
\CS_{s_i}(X,Y)=x_{i+1}+\cdots+x_n+y_1+\cdots+y_{n-i}.
\]

\medskip
\noin
4) Consider the partition $\l=\rho_k$ for some $k\lequ n$. Then
Theorem \ref{mgras} gives
\begin{equation}
\label{Dk}
\CS_{\rho_k}(X,Y)=\sum_{\a\in \D_k}\wt{Q}_{\a}(X)
\wt{Q}_{\rho_k\ssm\a}(Y).
\end{equation}
Further special cases of Theorem \ref{mgras} are given in Corollaries
\ref{ssc} and \ref{newcor}.

\subsection{Stability}
\label{stab}
For each $m<n$ let $i=i_{m,n}\colon  W_m \hra W_n$ be the
embedding via the first $m$ components. The maps $i_{m,n}$
are used to define the infinite hyperoctahedral group
$W_{\infty}=\bigcup_nW_n$. We say that a type $C$ double
Schubert polynomial $\CS_w(X_m,Y_m)$ with $w\in W_m$ is {\em stable} if
\[
\left.
\CS_{i(w)}(X_n,Y_n)\right|_{x_{m+1}=\cdots=x_n=y_{m+1}=\cdots =
y_n=0}\, =\,\CS_w(X_m,Y_m)
\]
for all $n>m$. Theorems \ref{mgras} and \ref{posit}
show that this stability
property fails in general, even for maximal Grassmannian elements.
For example when $m=2$, $n=3$ and $w=(\ov{2},1)$ we have
\[
\left. \CS_{i(w)}(X_3,Y_3)\right|_{x_3=y_3=0} =
x_1x_2+(x_1+x_2)(y_1+y_2)+y_1^2+y_1y_2+y_2^2
\]
which differs from $\CS_w(X_2,Y_2)$, given in Example 1. 
In addition, Theorem \ref{posit}
shows that the only permutation
$\om\in S_n$ for which $\CS_{\om}(X,Y)$ is stable is
$\om=1$.

For $u,v\in W_{\infty}$ we say that $u$ precedes $v$ in the
weak Bruhat order, and write $u \leq v$, if there are generators
$s_{a_1},\ldots,s_{a_r}\in W_{\infty}$ with $v = u s_{a_1}\cdots s_{a_r}$
and $r=\ell(v)-\ell(u)\gequ 0$. For each $k>0$ let $w_{\rho_k}$ be the
maximal Grassmannian permutation corresponding to $\rho_k$, that is
\[
w_{\rho_k}=[\ov{k}, \ov{k-1},\ldots,\ov{1},k+1,k+2,\ldots \, ]\in W_{\infty}.
\]

\begin{prop}[Stability]
\label{Cstab}
If $w\leq w_{\rho_k}$ for some $k>0$ then $\CS_w(X,Y)$ is stable.
\end{prop}
\begin{proof} 
Equation (\ref{Dk})
above shows that $\CS_w(X,Y)$ is stable when $w=w_{\rho_k}$. It follows
from the definition of the polynomials $\CS_w(X,Y)$ by divided difference
operators that
\[
\partial'_i\CS_w = \left\{ \begin{array}{cl}
             \CS_{ws_i} & \mathrm{ if } \ \ ws_i\leq w, \\
             0 &   \mathrm{ otherwise. }
       \end{array} \right.
\]
This implies that if $\CS_w(X,Y)$ is stable, then so is $\CS_u(X,Y)$ for
any $u\leq w$. \end{proof}

\medskip
The elements $w=(w_1,w_2,\ldots)\in W_{\infty}$ which satisfy the 
hypothesis of Proposition \ref{Cstab} are characterized by the following
three properties, which should hold for all $i<j$: (i) if $w_i$ and 
$w_j$ are unbarred then $w_i<w_j$, (ii) if
$w_i$ and $w_j$ are barred then $|w_i|>|w_j|$,
(iii) the barred entries are smaller than the unbarred entries
in absolute value.
One sees that the only maximal Grassmannian elements $w_{\l}$ 
with these properties are the $w_{\rho_k}$ themselves. 

Lascoux, 
Pragacz and Ratajski have shown in \cite[Thm.\ A.2]{LP1} that
the symplectic Schubert polynomials $\CS_w(X)$ are stable,
for any $w\in W_m$, in the sense that
\[
\left.
\CS_{i(w)}(X_n)\right|_{x_{m+1}=\cdots=x_n=0}\, =\, \CS_w(X_m).
\]
When $w$ is a maximal Grassmannian element $w_{\l}$, the stability
of $\CS_w(X)=\wt{Q}_{\l}(X)$ is clear. 
Stability for $\CS_{\om}(X)$ when $\om\in S_n$
is less obvious, and is equivalent (by Corollary \ref{CtoS})
to the following `shift' property
of type $A$ single Schubert polynomials, which can be checked directly.

\begin{cor}
\label{Ashift}
If $\om\in S_n$ with $\om(1)=1$ and $\tau$ is the $n$-cycle
$(n\dots 21)$, then
\[
\left.
\AS_{\om}(\tau X)\right|_{x_n=0}=\AS_{\tau\om\tau^{-1}}(X).
\]
\end{cor}

\subsection{Further properties}
The next result uses the Cauchy formula of Corollary \ref{Cauchy}
to give a characterization of the type $C$ double Schubert polynomials.

\begin{prop}
\label{char}
Suppose that we are given a homogeneous 
polynomial $P_{w_0}(X,Y)\in \Z[X,Y]$,
and for each $w\in W_n$ define $P_w(X,Y) 
= \partial'_{w^{-1}w_0} (P_{w_0}(X,Y))$. Let
$P_{\om}(X)=P_{\om}(X,0)$ and $P_{\l}(X)=
P_{w_{\l}}(X,0)$ for every pair $(\om,\l)\in S_n\times \D_n$.
Assume that these polynomials satisfy the following properties:

\medskip
{\em (i)} $P_{w_0}(X,Y) = \sum_{\om,\l}
P_{\om\om_0}(-\om_0 X) P_{\l'}(X) P_{\om}(\om_0 Y) P_{\l}(Y)$

\medskip
{\em (ii)} $P_{\om_0}(X)=\AS_{\om_0}(\om_0 X)=x_2x_3^2\cdots x_n^{n-1}$
\ and \ $P_{\rho_n}(X)=\wt{Q}_{\rho_n}(X)$.

\medskip
\noindent
Then $P_w(X,Y)=\CS_w(X,Y)$ for all $w\in W_n$.
\end{prop}
\begin{proof}
Recall that if a polynomial $f$ is homogeneous of degree $d$, then
$\partial_i'f$ is homogeneous of degree $d-1$ for each $i$. Hence, the
definition of the polynomials $P_w$ and their normalization in (ii)
imply that $P_w(X)$ is homogeneous of degree $\ell(w)$. In 
particular $P_w(0)=0$ unless $w=1$, when $P_1(X)=1$. 
Properties (i) (with $Y=0$) and (ii) thus give 
\[
P_{w_0}(X)=P_{\om_0}(-\om_0X)P_{\rho_n}(X)=(-1)^{n(n-1)/2}
\AS_{\om_0}(X)\wt{Q}_{\rho_n}(X)=\CS_{w_0}(X).
\]
We deduce that $P_w(X)=\partial'_{w^{-1}w_0} (P_{w_0}(X))=\CS_w(X)$
for all $w\in W_n$. The result now follows from (i) and Corollary
\ref{Cauchy}.
\end{proof}

\medskip

We close this section  with a vanishing property, which reflects the fact 
that the top polynomial $\CS_{w_0}(X,Y)$ represents the class of the
diagonal in flag bundles (see \cite[Thm. 1.1]{Gra}, \cite[Sect.\ 8]{Bi} 
and Section \ref{dl}).

\begin{prop}[Vanishing]
We have
\[
\CS_{w_0}(X,wX)=\left\{ \begin{array}{cl}
             2^nx_1\cdots x_n\prod_{i>j}(x_i^2-x_j^2)
             & \mathrm{ if } \ \ w=\om_0\in S_n, \\
             0 &   \mathrm{ otherwise. }
       \end{array} \right.
\]
\end{prop}
\begin{proof}
The vanishing property for $\wt{Q}(X,Y)$ from \cite[Sect.\ 2]{LP1} gives
\[
\wt{Q}(X,wX)=\left\{ \begin{array}{cl}
             \prod_{i\gequ j}(x_i+x_j) & \mathrm{ if } \ \ w\in S_n, \\
             0 &   \mathrm{ otherwise. }
       \end{array} \right.
\]
On the other hand, $\Delta(X,\om X)=0$ for all $\om\in S_n\ssm \{\om_0\}$.
\end{proof}

\section{Symplectic degeneracy loci}
\label{dl}

\subsection{Main representation theorem}
\label{firstss}
Let $V$ be a vector bundle of rank $2n$ on an algebraic variety $\X$.
Assume that $V$ is a {\em symplectic} bundle, i.e.\ $V$ is equipped 
with an everywhere nondegenerate skew-symmetric form $V\otimes V\ra
\C$. Consider two  Lagrangian (i.e., maximal isotropic) subbundles
$E$, $F$ of $V$ together with flags of subbundles
\begin{align*}
0=E_0 \subset E_1\subset E_2 \subset \cdots \subset
E_n =E\subset V  \\
0=F_0 \subset F_1\subset F_2 \subset \cdots \subset
F_n =F \subset V,
\end{align*}
where $\dim E_i = \dim F_i = i$ for all $i$. 
These can be extended to complete flags $E_{\bullet}$, $F_{\bullet}$
in $V$ by defining
$E_{n+i}=E_{n-i}^{\perp}$ for $1\lequ i\lequ n$, and similarly
for $F_{\bullet}$. 
For each $i$ with $1\lequ i \lequ n$ let
\[
x_i=-c_1(E_{n+1-i}/E_{n-i}) \ \ \ \ 
\mathrm{and} \ \ \ \
y_i=-c_1(F_i/F_{i-1}).
\]

There is a group monomorphism $\phi\colon W_n\hra S_{2n}$ with image
\[
\phi(W_n)=\{\,\s\in S_{2n} \ |\ \s(i)+\s(2n+1-i) = 2n+1,
 \ \ \mathrm{for} \ \mathrm{all} \ \ i\,\}.
\]
The map $\phi$ is determined by setting, 
for each $w=(w_1,\ldots,w_n)\in W_n$, 
\begin{equation}
\label{phidef}
\phi(w)(i)=\left\{ \begin{array}{cl}
             n+1-w_{n+1-i} & \mathrm{ if } \ w_{n+1-i} \ \mathrm{is} \ 
             \mathrm{unbarred}, \\
             n+\ov{w}_{n+1-i} & \mathrm{otherwise}.
             \end{array} \right.
\end{equation}
For every $w\in W_n$ define the {\em degeneracy
locus} $\X_w\subset \X$ as the locus of $a \in \X$ such that
\[
\dim(E_r(a)\cap F_s(a))\gequ \#\,\{\,i \lequ r \ |\ \phi(w)(i)> 2n-s\,\}
\ \ \mathrm{for} \ \ 1\lequ r\lequ n,\, 1\lequ s\lequ 2n.
\]
The precise definition of $\X_w$ as a subscheme of $\X$ is obtained 
by pulling back from the universal case, which takes place on the
bundle $F(V)\ra \X$ of isotropic flags in $V$, following \cite{F3}; 
see also \cite[Sect.\ 6.2 and App.\ A]{FP}. If $G_{\bullet}$ denotes
the universal flag over $F(V)$, then the 
flag $E_{\bullet}$ corresponds to a section $\s\colon \X\ra F(V)$ of 
$F(V)\ra \X$ such that $\s^*(G_{\bullet})=E_{\bullet}$. The 
degeneracy locus $\X_w\subset \X$ is defined as the inverse image
$\s^{-1}(\wt{\X}_w)$ of the {\em universal Schubert variety} $\wt{\X}_w$
in $F(V)$. 

Assume that 
$\X_w$ is of pure codimension $\ell(w)$ and $\X$ is 
Cohen-Macaulay. In this case, we have a formula for
the fundamental class $[\X_w]$ of $\X_w$ in the Chow
group $CH^{\ell(w)}(\X)$ of codimension $\ell(w)$ algebraic cycles
on $\X$ modulo rational equivalence, expressed as a 
polynomial in the Chern roots $X=\{x_i\}$ and $Y=\{y_i\}$ 
of the vector bundles
$E$ and $F$. The exact result is

\begin{thm}[Locus Representation]
\label{lagfulthm}
$[\X_w]=\CS_{w}(X,Y)$ in $CH^{\ell(w)}(\X)$.
\end{thm}

\begin{proof} We will show that the above statement is equivalent
to the corresponding result in \cite{F3} (Theorem 1). Note that,
as in loc.\ cit., the proof is obtained by pulling back the 
corresponding equality on the flag bundle $F(V)$. In order to
avoid notational confusion define new alphabets
$U=(u_1,\ldots,u_n)$ and $Z=(z_1,\ldots,z_n)$ with
\[
u_i=-c_1(E_i/E_{i-1})=x_{n+1-i} \ \ \ \ 
\mathrm{and} \ \ \ \
z_i=-c_1(F_{n+1-i}/F_{n-i})=y_{n+1-i},
\]
in agreement with the conventions in \cite{F3}. 
Identify the Weyl group $W_n$ with its image $\phi(W_n)\subset S_{2n}$
and let $\partial_1^u,\ldots,\partial_n^u$ denote
the divided difference operators with respect to the $u$ variables,
defined as in loc.\ cit. For example 
$\partial_n^u(f)=(f-s_nf)/(2u_n)$, where $s_n$ sends $u_n$ to $-u_n$.
Theorem 1 of \cite{F3} gives
\begin{equation}
\label{fthm}
[\X_w]=\partial^u_{\phi(w^{-1}w_0)} \left(
\Delta(U,Z)\,F(U,Z)\right),
\end{equation}
where 
\[
F(U,Z)=\det\left(e_{n+1+j-2i}(U)+e_{n+1+j-2i}(Z)\right)_{i,j}.
\]
Using the criterion in \cite[Theorem 1.1]{Gra}, we see that
one can replace the kernel $\Delta(U,Z)\,F(U,Z)$ in (\ref{fthm}) with
\begin{equation}
\label{grathm}
(-1)^{n(n-1)/2}\prod_{i>j}(u_i-y_j)\,
F(U,Y)=(-1)^{n(n-1)/2}\om_0(\Delta(U,Y))\,F(U,Y),
\end{equation}
where the permutation $\om_0$ acts on the $u$ variables
(compare with the analysis in \cite[Section 5]{Gra}). Next we
apply \cite[Proposition 2.4]{LP1} to replace $F(U,Y)$ 
in (\ref{grathm}) by $\wt{Q}(U,Y)$, as these two
polynomials are congruent modulo the relations in the Chow
ring of the flag bundle of $V$.
Note also that we have the identities
\begin{equation}
\label{rels}
\partial^u_{n-i}\om_0=-\om_0\partial^u_i, \ \ 1\lequ i\lequ n-1,
\ \ \ \ 
\mathrm{and} \ \ \ \ 
\partial^u_n\om_0=\om_0\partial^u_0,
\end{equation}
where we define $\partial^u_0(f)=(f-s_0f)/(2u_1)$. 
The map $\phi$ interchanges $s_i$ with $s_{n-i}$ for all $i$; 
it follows from (\ref{fthm}), (\ref{grathm}) and (\ref{rels}) that
\[
[\X_w]=(-1)^{n(n-1)/2}\om_0(\partial^u_{w^{-1}w_0})'\left(
\Delta(U,Y)\,\wt{Q}(U,Y)\right)=\CS_w(X,Y),
\]
as required.
\end{proof}

\subsection{Lagrangian degeneracy loci for subbundles}
We turn next to the special case of Theorem \ref{lagfulthm}
which describes Lagrangian degeneracy loci. Note that
the variables $X=\{x_1,\ldots,x_n\}$ are the Chern roots of $E^*$ and 
for each $i$, $Y_i=\{y_1,\ldots,y_i\}$ are the Chern roots of 
$F_i^*$. 
For each strict partition $\l\in \D_n$ define the {\em degeneracy
locus} $\X_{\l}\subset \X$ as
\begin{equation}
\label{Cdef}
\X_{\l}=\{\,a\in\X\ |\ \dim(E(a)\cap F_{n+1-\l_i}(a))
\gequ i \ \mathrm{for} \ 1\lequ i\lequ \ell(\l)\,\}.
\end{equation}
The precise scheme-theoretic definition $\X_{\l}$ is obtained as before, 
by pulling back the corresponding universal locus $\wt{\X}_{w_{\l}}$
by the section $\s$. There is an equality $\X_{\l}=\X_{w_{\l}}$ of schemes,
hence the analysis in Section \ref{firstss} applies.
Assuming that $\X_{\l}$ is of pure codimension $|\l|$ and $\X$ is 
Cohen-Macaulay, Theorem \ref{lagfulthm} says that the
fundamental class $[\X_{\l}]$ equals $\CS_{\l}(X,Y)$ in the Chow
group $CH^{|\l|}(\X)$. Now Theorem \ref{mgras} gives
\begin{cor}[Lagrangian Degeneracy Loci]
\label{lagcor}
We have
\[
[\X_{\l}] =(-1)^{e(\l)+|\l'|}
\sum_{\a}\wt{Q}_{\a}(E^*)
\sum_{\b}
(-1)^{e(\a,\b)+|\b|}
\wt{Q}_{(\a\cup \b)'}(F^*)
\det(c_{\b_i-\l'_j}(F^{*}_{n-\l'_j}))
\]
in $CH^{|\l|}(\X)$, 
where the first sum is over all $\a\in \D_n$ and the
second over $\b\in \D_n$ with $\b\supset \l'$,
$\ell({\b})=\ell({\l'})$ and $\a\cap\b=\emptyset$.
\end{cor}

We can apply Corollary \ref{lagcor} to
recover the results of \cite{PR} which compute
the class $[\X_{\l}]$ for some special partitions $\l$. For 
instance, when $\l=k$ for some $k\lequ n$ the locus $\X_{\l}$ is
given by a single Schubert condition. Let $s_r$ denote the Schur
polynomial (or complete
homogeneous function) corresponding to $r$. 

\begin{cor}[\cite{PR}, Proposition 6.1] We have
\label{ssc}
\[
[\X_k]=\sum_{p=0}^k c_p(E^*) s_{k-p}(F_{n+1-k}^*).
\]
\end{cor}
\begin{proof}
For each $p < k$ there is an equality of $(k-p)\times (k-p)$
matrices
\begin{equation}
\label{mprod}
\{e_{1+j-i}(Y_{n-k+j})\} = \{e_{1+j-i}(Y_{n-k+1})\}\cdot
\{e_{j-i}(Y_{n-k+j}\ssm Y_{n-k+1})\}.
\end{equation}
Taking determinants in (\ref{mprod}) shows that 
\begin{equation}
\label{apply}
\det(c_{1+j-i}(F^*_{n-k+j}))_{1\lequ i,j\lequ k-p}=
s_{k-p}(F^*_{n-k+1}).
\end{equation}
Now apply Corollary \ref{lagcor} and use (\ref{apply}) to obtain
\[
[\X_k]=\sum_{p=1}^k c_p(E^*) s_{k-p}(F_{n+1-k}^*)+
\sum_{q=1}^k(-1)^{q-1}c_q(F^*)s_{k-q}(F^*_{n+1-k}).
\]
The proof is completed by noting that
\[
\sum_{q=1}^k(-1)^{q-1}e_q(Y)s_{k-q}(Y_{n+1-k})=
\det(e_{1+j-i}(Y_{n+1-i}))_{1\lequ i,j\lequ k}=s_k(Y_{n+1-k}),
\]
which follows from (\ref{apply}) and the Laplace expansion of
a determinant.
\end{proof}

\medskip

\noin
When $\l=\rho_k$, Corollary \ref{lagcor} becomes Theorem 9.1 of
\cite{PR}. More generally, we have the following result:

\begin{cor}
\label{newcor}
Let $\l=\rho_k\ssm j$ for an integer $j$ with $1\lequ j \lequ k$. Then
\[
[\X_{\l}] = \sum_{p=j}^k(-1)^{k-p}c_{p-j}(F^*_{n-j})
\sum_{\substack{\a\in \D_k \\ \a\cap p=\emptyset}} (-1)^{e(\a,p)}
\wt{Q}_{\a}(E^*)\wt{Q}_{\rho_k\ssm(\a\cup p)}(F^*)
\]
where $e(\a,p)=\#\,\{\,i\ |\ \a_i>p\,\}$.
\end{cor}

\medskip
\noindent
{\bf Remark.} The formula for Lagrangian degeneracy loci
in Corollary \ref{lagcor} (and its predecessor from Theorem
\ref{mgras}) interpolates between two extreme cases, where it
assumes a simple form. These correspond to the partitions
$\l=\rho_k$ (Example 4 in Section \ref{A}) and $\l=k$ 
(Corollary \ref{ssc}), and were derived in \cite{PR} using purely 
geometric methods. 
Note that a very different formula
for the $\l=\rho_k$ locus is given in \cite{F2} \cite{F3}; it would be
interesting to have an analogous statement for the remaining maximal
Grassmannian loci.

\subsection{Restriction of type $A$ loci}
\label{inter}
We now give a geometric interpretation of Theorem \ref{posit}, in
the setting of flag bundles. Let $F'(V)$ be the partial $SL_{2n}$-flag
bundle of which parametrizes flags of subbundles 
\[
0=E'_0 \subset E'_1\subset E'_2 \subset \cdots \subset
E'_n \subset V
\]
with $\dim E'_i = i$ for all $i$;
by abuse of notation let $E'_{\bullet}$ also denote the universal flag 
of vector bundles over $F'(V)$. Suppose that 
\[
0=F'_0 \subset F'_1\subset F'_2 \subset \cdots \subset
F'_{2n}= V
\]
is a fixed complete flag of subbundles of $V$, and set
\[
x'_i=-c_1(E'_i/E'_{i-1}) \ \ \ \ 
\mathrm{and} \ \ \ \
y'_j=-c_1(F'_{2n+1-j}/F'_{2n-j})
\]
for $1\lequ i\lequ n$ and $1\lequ j\lequ 2n$. For any
permutation $\om\in S_n$, we have 
a universal codimension $\ell(\om)$ Schubert variety
$\X'_{\om}$ in
$F'(V)$, defined as the locus of $a\in F'(V)$ such that
\[
\dim(E'_r(a)\cap F'_s(a))\gequ \#\,\{\,i \lequ r \ |\ \om(i)> 2n-s\,\}
\ \ \mathrm{for} \ \ 1\lequ r\lequ n,\, 1\lequ s\lequ 2n.
\]
Identify $\om$ with a permutation in $S_{2n}$ by using the 
embedding $S_n\hra S_{2n}$ which fixes the last $n$ entries.
Then the class of $\X'_{\om}$ in $CH^{\ell(\om)}(F'(V))$ is given by
\[
[\X'_{\om}]=\AS_{\om}(X',Y')
\]
where $X'=(x_1',\ldots,x_n',0,\ldots,0)$ and $Y'=(y'_1,\ldots,y'_{2n})$.
This follows from \cite[Prop.\ 8.1]{F1} by dualizing; 
to navigate through the
changes in notation, we have found \cite[Sect.\ 2.2]{FP} helpful.

Let $\theta\colon F(V)\hra F'(V)$ be the natural inclusion and
$\theta^*\colon CH(F'(V))\ra CH(F(V))$ the induced homomorphism on
Chow groups. Define
the complete flag $F'_{\bullet}$ by setting $F'_j=F_j$ for 
$1\lequ j\lequ 2n$. We then have $\theta^*(x_i')=x_{n+1-i}$ for
$1 \lequ i \lequ n$, and hence can identify $\theta^*(X')$
with $\om_0X$, 
where $\om_0$ denotes the longest element in $S_n$. The symplectic form
on $V$ induces isomorphisms $V/(F_j^{\perp})\cong F_j^*$ for each $j$; it
follows that $\theta^*(y'_j)=-y_j$ for $1\lequ j\lequ n$. Using the
stability property of type $A$ double Schubert polynomials
\cite[(6.5)]{M1} and \cite[Cor.\ 2.11]{F1} we see that
\begin{equation}
\label{pullback}
\theta^*(\AS_{\om}(X',Y')) = \AS_{\om}(\om_0X,-Y)= \CS_{\om^*}(X,Y),
\end{equation}
where the last equality is Theorem \ref{posit}.

On the other hand, by comparing the
defining conditions for $\X'_{\om}$ and $\wt{\X}_{\om^*}$
one checks that 
\[
\theta^{-1}(\X'_{\om})= F(V)\cap \X'_{\om} = \wt{\X}_{\om^*}
\]
and that the intersection is proper and generically transverse along
$\wt{\X}_{\om^*}$. It follows that 
\begin{equation}
\label{restrict}
\theta^*[\X'_{\om}] = [\wt{\X}_{\om^*}]
\end{equation}
and thus Theorem \ref{posit} corresponds in geometry 
to the restriction of universal Schubert classes (\ref{restrict}).

\subsection{Schubert calculus}
\label{sc}
We now suppose that $V$ is a symplectic vector space and let
$\X=F_{Sp}(V)$ denote the variety
of complete isotropic flags in $V$.
There is a tautological 
isotropic flag $E_{\bullet}$ of vector bundles over $\X$, and we 
let $F_{\bullet}$ be a fixed isotropic
flag in $V$. For each $w\in W_n$ the locus $\X_w$ is a Schubert
variety in $\X$; 
let $\sigma_w=[\X_w]\in CH^{\ell(w)}(\X)$ be the corresponding
Schubert class. The classes 
$\{\sigma_w\}_{w\in W_n}$ form an integral 
 basis for the Chow ring $CH(\X)$. Hence,
for each $u,v,w\in W_n$, we can define {\em type $C$
structure constants} $f(u,v,w)$ 
by the equation
\begin{equation}
\label{Cstr}
\sigma_u\sigma_v=\sum_{w\in W_n} f(u,v,w) \sigma_w.
\end{equation}
The $f(u,v,w)$ are nonnegative integers.

Theorem \ref{lagfulthm} implies that $\s_w=\CS_w(X)$, that is, the 
symplectic Schubert polynomials represent the Schubert classes. We
deduce that the $\{\CS_w(X)\}_{w\in W_n}$ form a $\Z$-basis for the
quotient ring $\Z[X]/I_n$, where $I_n$ denotes the ideal of positive
degree Weyl group invariants in $\Z[X]$ (using the well known 
isomorphism of the latter with $CH(\X)$; see e.g.\ \cite{Bo}). 
Recall that 
the type $A$ Schubert polynomials $\AS_{\om}(X)$ for $\om\in S_n$ 
form a $\Z$-basis for the additive subgroup $H_n$ of $\Z[X]$ spanned by 
the monomials $x^{\a}$, $\a\subset\rho_{n-1}$ \cite[(4.11)]{M1}. 
Using Corollary \ref{CtoS} we obtain

\begin{cor}
Fix a permutation $\om\in S_n$.
Suppose that {\em (i)} the polynomial $P_{\om}\in \Z[X]$ represents
the Schubert class $\s_{\om}$ in $CH(\X) \cong \Z[X]/I_n$
and {\em (ii)} $P_{\om} \in \om_0 H_n$. Then 
$P_{\om}(X)=\CS_{\om}(X)=\AS_{\om^*}(\om_0 X)$.
\end{cor}

We will 
apply the results of Section \ref{A} to obtain some information about the
structure constants $f(u,v,w)$. For this we need the following

\begin{defn}
Let $u,v\in S_n$ be two permutations. We say that the pair $(u,v)$ is
{\em $n$-stable} if the product $\AS_u(X)\AS_v(X)$ is contained in $H_n$.
\end{defn}

\noindent
The previous remarks show that  $(u,v)$ is $n$-stable 
exactly when there is an equation of type $A$
Schubert polynomials
\begin{equation}
\label{Astr}
\AS_u(X)\AS_v(X)=\sum_{w\in S_n}c(u,v,w)\AS_w(X)
\end{equation}
for some integers $c(u,v,w)$ (which are type $A$ structure constants).

\medskip
\noin
{\bf Examples.} 1) If $u\in S_k$ and $v\in S_{\ell}$
then $(u,v)$ is $n$-stable for any $n\gequ k+\ell-1$.

\medskip
\noin
2) Let $\leq$ denote the weak Bruhat order. 
The Leibnitz rule for divided differences
\begin{equation}
\label{leibnitz}
\partial_i(fg)=(\partial_if)g+(s_if)(\partial_ig)
\end{equation}
implies that
if $(u,v)$ is $n$-stable and $u'\leq u$, $v'\leq v$
then $(u',v')$ is also $n$-stable.

\medskip
\noin
3) For each $i>0$ and $\om\in S_{\infty}$ let $d_i(\om)=
\deg_{\, x_i}(\AS_{\om}(X))$. Then $(u,v)$ is $n$-stable if and only if
$d_i(u)+d_i(v)\lequ \max\{n-i,0\}$ for all $i$. One sees that
$d_1(\om)$ equals the number of nonempty columns in the diagram
of $\om$, or alternatively
\[
d_1(\om)=\#\,\{\,j\ |\ \exists \  i<j \ :\ \om(i)>\om(j)\,\}.
\]
This follows e.g.\ from the combinatorial algorithm for constructing
Schubert polynomials due to Bergeron \cite[Chap.\ 4]{M1}.
In particular, $d_1(\om)\gequ \om(1)-1$ for all $\om$.
We deduce that if $u(1)+v(1)>n+1$ then the pair
$(u,v)$ is not $n$-stable. 

\medskip
\noindent
4) For any partition $\l$ (not necessarily strict) and integer $r\gequ 
\ell(\l)$ there is a Grassmannian permutation $\om_{\l,r}\in S_{\infty}$,
determined by $\om_{\l,r}(i)=\l_{r+1-i}+i$ for $i\lequ r$ and 
$\om_{\l,r}(i)<\om_{\l,r}(i+1)$ for $i>r$. One knows
\cite[(4.7)]{M2} that (for $n$ sufficiently large)
the Schubert polynomial
$\AS_{\om_{\l,r}}(X_n)$ is equal to the Schur function 
$s_{\l}(x_1,\ldots,x_r)$;
moreover $\deg_{x_i}(s_{\l})= \l_1$ for each $i$. Let
$\mu$ be a second partition and assume that the integer $s$ satisfies
$\ell(\mu)\lequ s \lequ r$. We deduce that the pair
$(\om_{\l,r}, \om_{\mu,s})$ is $n$-stable if and only if 
$n\gequ \l_1+\max\{r,\mu_1+s\}-1$.

\medskip
\noindent
5) If $\om\in S_n$ we denote by  $1\times \om$ the
permutation $(1,\om(1)+1,\ldots,\om(n)+1)$ in $S_{n+1}$. 
We then have the following

\begin{prop}
\label{nstable}
Suppose that $u,v\in S_n$. Then $(u,v)$ is $n$-stable if and only if
$(1\times u,1\times v)$ is $(n+1)$-stable.
\end{prop}
\begin{proof} Given the Schubert polynomial $\AS_{1\times u}(X)$
one can recover $\AS_u(X)$ by using
Corollary \ref{Ashift}. It follows easily from this 
that if $(1\times u,1\times v)$ is $(n+1)$-stable then
$(u,v)$ is $n$-stable. For the converse, apply
\cite[(4.22)]{M1} to obtain
\begin{equation}
\label{1timesw}
\AS_{1\times u}(X)=\partial_1\cdots \partial_{n-1}(x_1\cdots x_{n-1}\AS_u(X))
\end{equation}
for every $u\in S_n$. Assuming that $(u,v)$ is $n$-stable, we deduce that
$(1\times u,1\times v)$ is $(n+1)$-stable from (\ref{leibnitz}),
(\ref{1timesw}) and a straightforward induction.
\end{proof}

\medskip

Extend the involution $*\colon S_n\ra S_n$ 
of Section \ref{mts} to all of $W_n$ 
by letting $w^*=\om_0w\om_0$, for each $w\in W_n$. 

\begin{prop} 
\label{curprop}
Assume that $(u,v)\in S_n\times S_n$ is $n$-stable.
 Then for each $w\in W_n$,
\begin{equation}
\label{curious}
f(u^*,v^*,w^*)=\left\{ \begin{array}{cl}
             c(u,v,w) & \mathrm{ if } \ \ w\in S_n, \\
             0 &   \mathrm{ otherwise.}
       \end{array} \right.
\end{equation}
\end{prop}
\begin{proof}
We apply the permutation $\om_0$ to the polynomials in 
(\ref{Astr}), then use Corollary \ref{CtoS} to obtain
\[
\CS_{u^*}(X)\CS_{v^*}(X)=\sum_{w\in S_n}c(u,v,w)\CS_{w^*}(X),
\]
and compare this with (\ref{Cstr}).
\end{proof}

\medskip
\noindent
{\bf Remarks.} 1) If $u,v\in S_m$ then $(u,v)$ is $n$-stable 
for any $n\gequ 2m-1$, thus
(\ref{curious}) is true for all such $n$. 
Hence every type $A$ structure constant 
$c(u,v,w)$ appears as the stable limit of type $C$ constants
$f(u^*, v^*,w^*)$, with the involution $*$ defined using the
rank of the corresponding hyperoctahedral group.

\medskip
\noindent
2) Poincar\'e duality for the Schubert classes in the complete
$SL_n$-flag variety implies that
\[
c(u,v,w)=c(u, \om_0 w, \om_0v)=c(u^*,v^*,w^*)
\]
for all $u,v,w\in S_n$. Moreover, the structure constants $c(u,v,w)$
and $f(u,v,w)$ are stable under the inclusions of the respective
Weyl groups. These properties are compatible because of the general
identity
\begin{equation}
\label{ident}
c(1\times u, 1\times v, 1\times w)= c(u,v,w),
\end{equation}
valid for all $u,v,w\in S_{\infty}$.
Note that (\ref{ident}) is a special case of \cite[Cor.\ 4.5.6]{BS}.

\medskip
\noindent
3) The special case of the Pieri rule of \cite{PR1} displayed after
Example 2.3 of loc.\ cit.\  and the Giambelli-type formula in
\cite[Cor.\ 2.4]{PR1} are consequences of 
Proposition \ref{curprop}, applied to (non-maximal) isotropic
Grassmannians. To recover these
formulas from 
Proposition \ref{curprop}, notice that if $\om_{\l}\in S_n$ is
a Grassmannian permutation corresponding to the partition $\l$,
then $\om_0\om_{\l}\om_0$ is also Grassmannian and
corresponds to the conjugate (or transpose) of $\l$.

\vspace{0.2cm}

\begin{table}
\begin{tabular}{|c|c|c|} \hline

&& \\
Fibration & $F_{Sp}(2n)\to LG(n,2n)=LG$ & $F_{SL}(n) \times LG \to LG$ \\

&&\\

Basis & $\{\CS_w(X)\}_{w\in W_n}$ & 
$\{\CS_{\om,\l}(X)\}_{(\om,\l)\in S_n \times \D_n}$ \\

&& \\
Group & $S_n \ltimes (\Z_2)^n=W_n$ & $S_n \times (\Z_2)^n$ \\ 

&&\\ \hline

\end{tabular} 
\skipline
\caption{}
\label{table}
\end{table}

Table \ref{table} 
illustrates the connections between 
some of the geometric, algebraic and combinatorial
objects in this paper. There is 
a fibration of the symplectic flag variety $F_{Sp}(2n)=Sp(2n)/B_1$ over
the Lagrangian Grassmannian $LG(n,2n)=Sp(2n)/P_n$, 
with fiber equal to the $SL_n$-flag
variety $F_{SL}(n)=SL(n)/B_2$. The right column in Table \ref{table}
corresponds to the trivial product fibration 
$F_{SL}(n)\times LG(n,2n)$. The symplectic Schubert polynomials
give a $\Z$-basis for the cohomology (or Chow) ring in either column, 
as shown. When restricted to the parameter
spaces $S_n$ and $\D_n$, the Schubert polynomial and product 
bases coincide.
The correspondence between the set $\D_n$ 
and the group $(\Z_2)^n=\{\pm 1\}^n$ is given as follows:
the strict partition $\l\in \D_n$ corresponds to the vector
$(\epsilon_r)\in\{\pm 1\}^n$ with
\[
\epsilon_r= \left\{ \begin{array}{cl}
           -1 
           & \mathrm{ if } \ r \ \mathrm{ is } \ \mathrm{a}
             \ \mathrm{ part } \ \mathrm{of} \ \lambda, \\
           1 
           & \mathrm{ otherwise. } 
             \end{array} \right.
\]

\section{Orthogonal Schubert polynomials and degeneracy loci}
\label{C}

\noindent
In this section we describe the analogues of the main results of 
Sections \ref{A} and \ref{dl} for the orthogonal groups
of types $B_n$ and $D_n$. The alphabets $X_k$ and $Y_k$ for
$1 \lequ k\lequ n$ are 
defined as before, with $X=X_n$ and $Y=Y_n$. For each
$\l\in \D_n$ define the $\wt{P}$-polynomial
\[
\wt{P}_{\l}(X)=2^{-\ell(\l)}\wt{Q}_{\l}(X)
\]
and, for $k\in\{n-1,n\}$, the reproducing kernel
\[
\wt{P}_k(X,Y)=
\sum_{\l\in \D_k}\wt{P}_{\l}(X)\wt{P}_{\rho_k\ssm\l}(Y).
\]

\subsection{Type $B$ double Schubert polynomials and loci}
\label{oddorth}

The Weyl group of the root system $B_n$ 
is the same as that for $C_n$, as are the
divided difference operators $\partial_i$ for 
$1\lequ i \lequ n-1$. However the operator $\partial_0$ differs
from the one in type $C$ by a factor of $2$:
\[
\partial_0(f)=(f-s_0f)/x_1
\]
for any $f\in A[X]$. In this and the next subsection we will 
take $A=\Z[1/2][Y]$.
The $\partial_i$ are used as before
to define operators 
$\partial'_w$ for all barred permutations $w\in W_n$. 

\begin{defn}
For every $w\in W_n$ the {\em type $B$ 
double Schubert polynomial $\BS_w(X,Y)$} is given by
\[
\BS_w(X,Y)= (-1)^{n(n-1)/2}\partial'_{w^{-1}w_0}\left(
\Delta(X,Y)\wt{P}_n(X,Y)\right).
\]
\end{defn}

For each partition $\l\in \D_n$ let
$\BS_{\l}(X,Y)=\BS_{w_{\l}}(X,Y)$.
Using the results of \cite[Appendix]{LP2} and the same arguments
as in Section \ref{A} produces the following two theorems.

\begin{thm}[Maximal Grassmannian]
\label{mgrasB}
For any strict partition
$\l\in \D_n$ the double Schubert polynomial
$\BS_{\l}(X,Y)$ is equal to 
\[
(-1)^{e(\l)+|\l'|}
\sum_{\a}\wt{P}_{\a}(X)
\sum_{\b}
(-1)^{e(\a,\b)+|\b|}
\wt{P}_{(\a\cup \b)'}(Y)
\det(e_{\b_i-\l'_j}(Y_{n-\l'_j})),
\]
where the first sum is over all $\a\in \D_n$ and the
second over $\b\in \D_n$ with $\b\supset \l'$,
$\ell({\b})=\ell({\l'})$ and $\a\cap\b=\emptyset$.
\end{thm}

\begin{thm}[Positivity]
\label{positB}
For every $\om\in S_n$,
\[
\BS_{\om}(X,Y)=\AS_{\om^*}(\om_0X,-Y).
\]
In particular, $\BS_{\om}(X,Y)$ has nonnegative integer coefficients.
\end{thm}

\medskip
Consider an orthogonal vector bundle $V$ of rank $2n+1$ 
on an algebraic variety $\X$, so that $V$ is 
equipped with an everywhere nondegenerate quadratic form. 
One is given flags
of isotropic subbundles of $V$ as before, with $E=E_n$ and $F=F_n$. 
These are extended to complete flags $E_{\bullet}$ and
$F_{\bullet}$ by setting
$E_{n+i}=E_{n+1-i}^{\perp}$ and $F_{n+i}=F_{n+1-i}^{\perp}$
for $1\lequ i\lequ n+1$. The classes $x_i,y_i\in CH^1(\X)$
are defined as in Section \ref{dl}. 
In this case we have a monomorphism $\psi\colon W_n\hra S_{2n+1}$ with image
\[
\psi(W_n)=\{\,\s\in S_{2n+1} \ |\ \s(i)+\s(2n+2-i) = 2n+2,
 \ \ \mathrm{for} \ \mathrm{all} \ \ i\,\},
\]
determined by the equalities
\[
\psi(w)(i)=\left\{ \begin{array}{cl}
             n+1-w_{n+1-i} & \mathrm{ if } \ w_{n+1-i} \ \mathrm{is} \ 
             \mathrm{unbarred}, \\
             n+1+\ov{w}_{n+1-i} & \mathrm{otherwise}
             \end{array} \right.
\]
for each $w=(w_1,\ldots,w_n)\in W_n$.
The {\em degeneracy
locus} $\X_w$ is the locus of $a \in \X$ such that
\[
\dim(E_r(a)\cap F_s(a))\gequ \#\,\{\,i \lequ r \ |\ \psi(w)(i)> 2n+1-s\,\}
\]
for $1\lequ r\lequ n$, $1\lequ s\lequ 2n$.
Assuming that $\X_w$ is of pure codimension $\ell(w)$ and $\X$ is 
Cohen-Macaulay, we have

\begin{thm}[Locus Representation]
\label{fulthmB}
$[\X_w]=\BS_{w}(X,Y)$ in $CH^{\ell(w)}(\X)$.
\end{thm}
\noin
The proof of this theorem is the same as its analogue in Section
\ref{dl}, using the corresponding results of \cite{F3}, \cite{Gra}
and \cite{LP2}.
The maximal isotropic degeneracy loci $\X_{\l}$ for $\l\in \D_n$ are
defined by the same inequalities (\ref{Cdef}) as in the 
symplectic case. By combining the two previous theorems we
immediately obtain a formula for the class $[\X_{\l}]$ in
$CH^{|\l|}(\X)$:

\begin{cor}
\label{lagcorB}
We have
\[
[\X_{\l}] =(-1)^{e(\l)+|\l'|}
\sum_{\a}\wt{P}_{\a}(E^*)
\sum_{\b}
(-1)^{e(\a,\b)+|\b|}
\wt{P}_{(\a\cup \b)'}(F^*)
\det(c_{\b_i-\l'_j}(F^{*}_{n-\l'_j})),
\]
where the first sum is over all $\a\in \D_n$ and the
second over $\b\in \D_n$ with $\b\supset \l'$,
$\ell({\b})=\ell({\l'})$ and $\a\cap\b=\emptyset$.
\end{cor}

\subsection{Type $D$ double Schubert polynomials and loci}
\label{evenorth}

The situation here differs significantly from that of the previous
sections, so we will give more details.
The Weyl group $\wt{W}_n$ of type $D$ is an extension of $S_n$ by
an element $s_{\Box}$ which acts on the right by 
\[
(u_1,u_2,\ldots,u_n)s_{\Box}=(\ov{u}_2,\ov{u}_1,u_3,\ldots,u_n).
\]
$\wt{W}_n$ may be realized as a subgroup of $W_n$ by sending 
$s_{\Box}$ to $s_0s_1s_0$. The barred permutation
$\wt{w}_0\in\wt{W}_n$ of maximal length is given by
\[
\wt{w}_0=\left\{ \begin{array}{cl}
           (\ov{1},\ldots,\ov{n}) & \mathrm{ if } \ n \ \mathrm{is} \ 
             \mathrm{even}, \\
           (1,\ov{2},\ldots,\ov{n}) & \mathrm{ if } \ n \ \mathrm{is} \ 
             \mathrm{odd}.
             \end{array} \right.
\]

For $1\lequ i \lequ n-1$ the
action of the generators $s_i$ 
on the polynomial ring $A[X]$ and the
divided difference operators $\partial_i$ 
are the same as before. We let $s_{\Box}$
act by sending $(x_1,x_2)$ to $(-x_2,-x_1)$ and fixing the
remaining variables, while 
\[
\partial_{\Box}(f):=(f-s_{\Box}f)/(x_1+x_2)
\]
for all $f\in A[X]$. Define operators $\partial'_i=-\partial_i$
for each $i$ with $1\lequ i\lequ n-1$ and set 
$\partial'_{\Box}:=\partial_{\Box}$; these are used as above to
define $\partial_w$ and $\partial'_w$ for all $w\in \wt{W}_n$.

\begin{defn}
For every $w\in \wt{W}_n$ the {\em type $D$ 
double Schubert polynomial $\DS_w(X,Y)$} is given by
\[
\DS_w(X,Y)= (-1)^{n(n-1)/2}\partial'_{w^{-1}\wt{w}_0}\left(
\Delta(X,Y)\wt{P}_{n-1}(X,Y)\right).
\]
\end{defn}

The maximal Grassmannian elements $w_{\l}$ 
in $\wt{W}_n$ are parametrized by partitions $\l\in \D_{n-1}$. 
For each such $\l$ we set $\ell=\ell(\l)$; then
\[
w_{\l}=(\ov{\l_1+1},\ldots,\ov{\l_{\ell}+1},\wh{1},\m_{n-\ell-1},
\ldots,\m_1)
\]
where $\m=\rho_n\ssm(\l_1+1,\ldots,\l_\ell+1,1)$ and 
$\wh{1}$ is equal to $1$ or $\ov{1}$ according to the parity
of $\ell$. Let 
$\l'=\rho_{n-1}\ssm \l$ be the dual partition of $\l$ and 
\[
k=\left\{ \begin{array}{cl}
           n-\ell & \mathrm{ if } \,\ n = \ell \ (\mathrm{mod}\ 2), \\
           n-\ell-1 & \mathrm{ if } \,\ n \neq \ell \ (\mathrm{mod}\ 2).
             \end{array} \right.
\]

Let $\D_{n-1}^+$ be the set of strictly decreasing sequences $\b$
of elements of the set $\{0,1,\ldots,n-1\}$, with length $\ell(\b)$
equal to
the number of terms in the sequence. We think of $\D_{n-1}^+$ as
consisting of partitions in $\D_{n-1}$, possibly with an `extra zero
part'. If $\a\in \D_{n-1}$ and $\b\in \D_{n-1}^+$ let
\[
n_i(\a,\b)=\#\,\{\,j\ |\ \a_i>\b_j\gequ\a_{i+1} \,\}
\ \ \ \ \mathrm{and} \ \ \ \ 
f(\a,\b)=\sum_{i=1}^{\ell(\a)} i\,n_i(\a,\b).
\]
Also define
\[
f(\l)=\left\{ \begin{array}{cl}
           e(\l)+\ell(\l)
           & \mathrm{ if } \,\ n = \ell \ (\mathrm{mod}\ 2), \\
           e(\l)  & \mathrm{ if } \,\ n \neq \ell \ (\mathrm{mod}\ 2)
             \end{array} \right.
\]
and note that $f(\l)=f(\l,\l')$, provided that $\l'$ is identified with an
element in $\D_{n-1}^+$ {\em of length} $k$. Let 
$\DS_{\l}(X,Y)=\DS_{w_{\l}}(X,Y)$; we can now state

\begin{thm}[Maximal Grassmannian I]
\label{mgrasD}
For any strict partition
$\l\in \D_{n-1}$ the double Schubert polynomial
$\DS_{\l}(X,Y)$ is equal to 
\[
(-1)^{f(\l)+|\l'|}
\sum_{\a}\wt{P}_{\a}(X)
\sum_{\b}
(-1)^{f(\a,\b)+|\b|}
\wt{P}_{(\a\cup \b)'}(Y)
\det(e_{\b_i-\l'_j}(Y_{n-1-\l'_j})),
\]
where the first sum is over all $\a\in \D_{n-1}$ and the
second over $\b\in \D^+_{n-1}$ with $\b\supset \l'$,
$\ell({\b})=k$ and $\a\cap\b=\emptyset$.
\end{thm}
\begin{proof}
Given $\l\in \D_{n-1}$ define $\ell$ and $k$ as above and 
associate to $\l$ a partition $\nu$ of length $k$
as follows:
\[
\n=\left\{ \begin{array}{cl}
           \rho_n\ssm (\l_1+1,\ldots,\l_\ell+1)
           & \mathrm{ if } \,\ n = \ell \ (\mathrm{mod}\ 2), \\
           \rho_n\ssm (\l_1+1,\ldots,\l_\ell+1,1)
           & \mathrm{ if } \,\ n \neq \ell \ (\mathrm{mod}\ 2).
             \end{array} \right.
\]
Note that $w_{\l}=\wt{w}_0\tau_{\l}\,\om_k\,\delta_k^{-1}$, where 
\begin{gather*}
\tau_{\l} = \left\{ \begin{array}{cl}
           (\n_k,\ldots,\n_1,\l_1+1,\ldots,\l_{\ell}+1)
           & \mathrm{ if } \,\ n = \ell \ (\mathrm{mod}\ 2), \\
           (\n_k,\ldots,\n_1,\l_1+1,\ldots,\l_{\ell}+1,1)
           & \mathrm{ if } \,\ n \neq \ell \ (\mathrm{mod}\ 2)
             \end{array} \right. \\
\om_k = (k,\ldots ,2,1,k+1,\ldots ,n) \\
\delta_k = (s_{n-k}\cdots s_2s_1s_{n-k+1}\cdots s_2s_{\Box})
           \cdots (s_{n-2}\cdots s_2s_1s_{n-1}\cdots s_2s_{\Box});
\end{gather*}
hence
\[
\partial_{w_{\l}^{-1}\wt{w}_0}=\partial_{\delta_k}\circ
\partial_{\om_k}\circ\partial_{\tau_{\l}^{-1}}.
\]
We compute as in the proof of Theorem \ref{mgras}:
\begin{align}
\partial_{\tau_{\l}^{-1}}(\Delta(X,Y))
&= \prod_{j=1}^k\prod_{p=1}^{n-\n_j}(x_{k+1-j}-y_p) \\
\label{e4D}
&= \sum_{\gamma}(-1)^{|\gamma|-|\n|}
\prod_{j=1}^k x_{k+1-j}^{n-\gamma_j}e_{\gamma_j-\n_j}
(Y_{n-\n_j})
\end{align}
where the sum (\ref{e4D}) is over all $k$-tuples $\gamma=(\gamma_1,
\ldots,\gamma_k)$ of nonnegative integers.

 Observe that
\[
\partial'_{\delta_k}\partial'_{\om_k}=
\partial_{\delta_k}\partial'_{\om_k}=
(-1)^{k/2}\partial_{\delta_k}\partial_{\om_k}
\]
since $k(k-1)/2 \equiv k/2$ (mod $2)$. The operator
$\partial_{\om_k}$ is a Jacobi symmetrizer; hence for any partition
$\theta=(\theta_i)$ we have
$\partial_{\om_k}(x_1^{\theta_1}\cdots x_k^{\theta_k})=0$ unless
$\theta_1>\theta_2>\cdots >\theta_k$, when
\[
\partial_{\om_k}(x_1^{\theta_1}\cdots x_k^{\theta_k})=
s_{\theta_1-k+1,\theta_2-k+2,\ldots ,\theta_k}(x_1,\ldots,x_k)
\]
is a Schur $S$-polynomial. In the latter case we can
apply \cite[Theorem 11]{LP2} and deduce that
for any partition $\l\in \D_{n-1}$,
\[
(-1)^{k/2}
\partial_{\delta_k}\partial_{\om_k}(x_1^{\theta_1}\cdots x_k^{\theta_k}
\wt{P}_{\l}(X))=0
\]
unless each part of $\ov{\theta}=
(n-1-\theta_k,\ldots,n-1-\theta_1)\in \D_{n-1}^+$ 
occurs in $\l$. In this case, the image is 
$(-1)^{f(\a,\ov{\theta})}\wt{P}_{\a}(X)$ where $\alpha=\l\ssm\ov{\theta}$
(note that the operator $\partial_{\Box}$ in \cite{LP2} differs from
ours by a sign). Moreover, if
$(m_1,\ldots,m_k)$ is a $k$-tuple of distinct nonnegative integers
and $\s\in S_k$ is such that $m_{\s(1)}>\cdots >m_{\s(k)}$ then
\[
\partial_{\om_k}(x_1^{m_1}\cdots x_k^{m_k})=
\mathrm{sgn}(\s)\partial_{\om_k}(x_1^{m_{\s(1)}}\cdots x_k^{m_{\s(k)}}),
\]
hence the previous analysis applies, up to $\mathrm{sgn}(\s)$. 

Noting that $\n_j-1=\l'_j$ for $1\lequ j\lequ k$ and 
$|\n|\equiv |\l'|\,($mod $2)$, 
use (\ref{e4D}) to compute
\begin{gather*}
\partial'_{w_{\l}^{-1}\wt{w}_0}(\Delta(X,Y)\, \wt{P}_{n-1}(X,Y))=
\partial'_{\delta_k}\partial'_{\om_k}\left(\wt{P}_{n-1}(X,Y) 
\partial'_{\tau_{\l}^{-1}}(\Delta(X,Y))\right) \\
= \partial'_{\delta_k}\partial'_{\om_k}
\Bigg(\wt{P}_{n-1}(X,Y) 
\sum_{\substack{\gamma\in \D_n  \\  \ell(\gamma)=k}}
(-1)^{|\gamma|-|\n|+\ell(\tau_{\l})}
\sum_{\s\in S_k}\prod_{j=1}^k x_{k+1-j}^{n-\gamma_{\s(j)}}
e_{\gamma_{\s(j)}-\n_j} (Y_{n-\n_j})\Bigg) \\
= (-1)^{r(\l)}\sum_{\a}\wt{P}_{\a}(X)
\sum_{\b}(-1)^{f(\a,\b)+|\b|}
\wt{P}_{(\a\cup \b)'}(Y)
\det(e_{\b_i-\l'_j}(Y_{n-1-\l'_j}))
\end{gather*}
where the ranges of summation are as in the statement of the
theorem and $r(\l)=\ell(\tau_{\l})+k(k-1)/2+|\l'|$. Finally, observe that
\[
\ell(\tau_{\l})=|\n|+\binom{\ell}{2}-
\binom{k+1}{2}\equiv f(\l)+\binom{k}{2}+\binom{n}{2} 
\ (\mathrm{mod} \ 2)
\]
so the signs fit to complete the proof. 
\end{proof}

\medskip
The analogue of Theorem \ref{posit} in type $D_n$ differs from the
one in type $B_n$. Let $H\cong(\Z_2)^{n-1}$ be the normal subgroup of
$\wt{W}_n$ consisting of those
elements equal to $(1,\ldots,n)$ in absolute
value. For each $\om\in S_n$, define the parameter space
\[
L(\om)=\{w\in \om H\ |\ \ell(w)=\ell(\om)\}.
\]

\begin{thm}[Positivity]
\label{positD}
For every $\om\in S_n$,
\begin{equation}
\label{Dsum}
\sum_{w\in L(\om)}\DS_w(X,Y)=\AS_{\om^*}(\om_0X,-Y).
\end{equation}
\end{thm}
\begin{proof}
We argue along the lines of the proof of Theorem \ref{posit},
using the operator 
\[
\Pi=(\partial_{\Box}+\partial'_1)\,\partial'_2\,
       (\partial_{\Box}+\partial'_1)\,(\partial'_3\partial'_2)\,
       \cdots\,
       (\partial_{\Box}+\partial'_1)\,(\partial'_{n-1}\cdots\partial'_2)\,
       (\partial_{\Box}+\partial'_1)
\]
in place of $\partial'_{v_0}$. We claim that
the image of $\Pi\colon \Z[1/2][X]\ra\Z[1/2][X]$ is a free 
$\Z[1/2][X]^{\wt{W}_n}$-module with basis $\{\AS_{\om}
(\om_0X)\}_{\om\in S_n}$. This statement is a type $D$ analogue of 
\cite[Prop.\ 4.1]{LP1}, and the proof is similar. The only 
difference is that here we use the relation 
$(\partial_{\Box}-\partial_1')\,\Pi =0$, and thus 
$\partial'_{w_{\l}s_{\Box}}(\partial_{\Box}-\partial_1')\,\Pi =0$,
for each $\l\in\D_{n-1}$. Moreover, \cite[Cor.\ 16]{LP2} implies that
\[
\partial'_{w_{\l}s_{\Box}}(\partial_{\Box}-\partial_1')(\wt{P}_{\m}(X))=
\left\{ \begin{array}{cl}
           1
           & \mathrm{ if } \,\ \m=\l, \\
           0
           & \mathrm{ if } \,\ |\m|=|\l| \ \, \mathrm{but} \ \, \m\neq\l
             \end{array} \right.
\]
(note that $\partial'_{w_{\l}s_{\Box}}\partial_{\Box}=
\partial'_{w_{\l}}$).
Now the vanishing property for $\wt{P}_{n-1}(X,Y)$ from \cite[Prop.\ 2]{LP2}
and the same argument as in \cite[Sect.\ 4]{LP1} show that
\[
\Pi(\AS_u(X)\wt{P}_{n-1}(X,Y))=
(-1)^{\ell(u)}\AS_u(\om_0 X)
\]
for every $u\in S_n$. 

Observe, since $L(\om_0)=\om_0H$, 
 that $\Pi=\sum_{w\in L(\om_0)}\partial'_{w^{-1}w_0}$.
We deduce as in the proof of Theorem \ref{posit}
that
\begin{equation}
\label{middle}
\sum_{w\in L(\om_0)}\DS_w(X,Y)=
(-1)^{\ell(\om_0)}\Pi(\Delta(X,Y)\wt{P}_{n-1}(X,Y))=
\AS_{\om_0}(\om_0X,-Y).
\end{equation}
For any permutation $\om$, apply the operator $\partial'_{\om^{-1}\om_0}$
to both sides of equation (\ref{middle}). We have seen that
\[
\partial'_{\om^{-1}\om_0}(\AS_{\om_0}(\om_0X,-Y))=
\AS_{\om^*}(\om_0X,-Y),
\]
so it remains to show that
\[
\partial'_{\om^{-1}\om_0}\sum_{w\in L(\om_0)}\DS_w(X,Y)=
\sum_{w\in L(\om)}\DS_w(X,Y)
\]
for all $\om\in S_n$. To prove this, note that for each $w\in L(\om_0)$,
the operator
$\partial'_{\om^{-1}\om_0}\partial'_{w^{-1}w_0}$ either vanishes or
equals $\partial'_{u^{-1}w_0}$, for a unique $u\in L(\om)$.
Moreover, different elements $w$ lead to different elements $u$.
\end{proof}

\smallskip

Note that the number of terms in the sum (\ref{Dsum}) equals $2^{h(\om)}$,
where $h(\om)$ is defined to be
the number of $j\gequ 2$ such that $\om(i)>\om(j)$ for all $i<j$, that
is, the number of `new lows' in the sequence $\om(1),\ldots,\om(n)$.

\begin{cor}
If $\om\in S_n$ satisfies $\om(1)=1$, then
\[
\DS_{\om}(X,Y)=\AS_{\om^*}(\om_0X,-Y).
\]
In particular, $\DS_{\om}(X,Y)$ has nonnegative integer coefficients.
\end{cor}

Consider now a vector bundle $V$ of rank $2n$
on an algebraic variety $\X$ with a quadratic form and
rank $n$ isotropic subbundles $E$ and $F$ with complete
flags of subbundles as in the type $C$ setting.
We assume that $E$ and $F$ are in the same family, 
that is $\dim(E(a)\cap F(a))\equiv n$ (mod $2$) for
every $a\in \X$. Define the classes 
$x_i,y_i\in CH^1(\X)$ as before.

There is a monomorphism $\phi\colon \wt{W}_n\hra S_{2n}$ whose
image consists of those permutations $\s\in S_{2n}$
such that $\s(i)+\s(2n+1-i) = 2n+1$ for all $i$ and the
number of $i\lequ n$ such that $\s(i)>n$ is even. The map
$\phi$ is defined by the same equation (\ref{phidef}) as in
the type $C$ case. 
Set $\delta(w)=0$ if $\ov{1}$ is a part of $w$, and
$\delta(w)=1$ otherwise, and define $\wt{\phi}\colon \wt{W}_n\hra S_{2n}$ 
by $$\wt{\phi}(w)=s_n^{\delta(w)}\phi(w).$$
The map $\wt{\phi}$ is a modification of 
$\phi$ so that in the sequence of values of $\wt{\phi}(w)$, 
$n+1$ always comes before $n$.
We need also the alternate complete flag $\wt{F}_{\bullet}$, with
$\wt{F_i}=F_i$ for $i\lequ n-1$ but completed with a 
maximal isotropic
subbundle $\wt{F}_n$ in the opposite family from $E$.  Define
\[
F_{\bullet}^{\delta}=\left\{ \begin{array}{cl}
           F_{\bullet}   & \mathrm{ if } \,\ n = \delta \ (\mathrm{mod}\ 2), \\
           \wt{F}_{\bullet}  & \mathrm{ if } \,\ n \neq 
           \delta \ (\mathrm{mod}\ 2).
             \end{array} \right.
\]

For $w\in \wt{W}_n$,
the {\em degeneracy locus} $\X_w$
is the locus of $a \in \X$ such that
\[
\dim(E_r(a)\cap F^{\delta(w)}_s(a))\gequ 
\#\,\{\, i \lequ r \ |\ \wt{\phi}(w)(i)>
2n-s\,\}
\]
for $1\lequ r\lequ n$, $1\lequ s\lequ 2n$. 
Recall that the flag $E_{\bullet}$ corresponds to a section 
$\sigma\colon  \X\ra F(V)$ of the bundle $F(V)\ra\X$ of isotropic flags in $V$. 
The subscheme $\X_w$ of $\X$ is then the inverse image under $\sigma$
of the closure 
of the locus of $y\in F(V)$ such that
\[
\dim(E_r(y)\cap F_s(y))= \#\,\{\, i \lequ r \ |\ \phi(w_0ww_0)(i)> 2n-s\,\}
\]
for $1\lequ r\lequ n-1$, $1\lequ s\lequ 2n$. This relates the present
formalism to that in \cite{F3}.
With the same assumptions on $\X$ and the codimension of 
$\X_w$ as before, and with the same arguments, we obtain

\begin{thm}[Locus Representation]
\label{fulthmD}
$[\X_w]=\DS_{w}(X,Y)$ in $CH^{\ell(w)}(\X)$.
\end{thm}

The maximal isotropic degeneracy locus $\X_{\l}$ for $\l\in \D_{n-1}$ 
is defined as
\begin{equation}
\label{maxlocD}
\X_{\l}= \{\,a\in\X\ |\ 
\dim(E(a)\cap F_{n-\l_i}(a))\gequ i \ \ \ \mathrm{for} \ \ \
1\lequ i\lequ \ell(\l)\,\}.
\end{equation}
Theorems \ref{mgrasD} and \ref{fulthmD} imply the following
formula for the class $[\X_{\l}]$ in $CH^{|\l|}(\X)$:

\begin{cor}
\label{lagcorD}
We have
\[
[\X_{\l}] =(-1)^{f(\l)+|\l'|}
\sum_{\a}\wt{P}_{\a}(E^*)
\sum_{\b}
(-1)^{f(\a,\b)+|\b|}
\wt{P}_{(\a\cup \b)'}(F^*)
\det(c_{\b_i-\l'_j}(F^{*}_{n-1-\l'_j})),
\]
where the first sum is over all $\a\in \D_{n-1}$ and the
second over $\b\in \D^+_{n-1}$ with $\b\supset \l'$,
$\ell({\b})=k$ and $\a\cap\b=\emptyset$.
\end{cor}

\subsection{Further results} The Propositions, Corollaries and Examples
in Sections \ref{A} and
\ref{dl} have orthogonal analogues. In particular there are
geometric interpretations of Theorems \ref{positB} and \ref{positD} 
(as in Section \ref{inter}) and connections
to the formulas of \cite{PR} for types $B$ and $D$. We omit most of
them here because their statements and proofs are straightforward, 
following the type $C$ case. The specialization
$Y=0$ produces type $B$ and $D$ {\em Schubert polynomials}
\begin{equation}
\label{BD}
\BS_w(X)=\BS_w(X,0) \ \ \ \ \ \ \mathrm{and} \ \ \ \ \ 
\DS_w(X)=\DS_w(X,0).
\end{equation}
Note that these polynomials differ from the orthogonal Schubert
polynomials defined in \cite{LP2} by a sign, which depends on the
degree. It is however still
true (arguing in the same way as \cite[Thm.\ A.2]{LP1})
that the polynomials (\ref{BD}) have a stability 
property: for any $w\in W_m$,
\[
\left.
\BS_{i(w)}(X_n)\right|_{x_{m+1}=\cdots=x_n=0}\,=\,\BS_w(X_m)
\]
where $i\colon  W_m \hra W_n$ is the natural embedding. In addition, the
set $\{\BS_{\om}(X)\BS_{\l}(X)\}$ for $\om\in S_n$, 
$\l\in\D_n$ forms an orthogonal product basis
of $A[X]$ with respect to the 
$A[X]^{W_n}$-linear scalar product
\[
A[X]\times A[X] \lra A[X]^{W_n}
\]
defined by the maximal divided difference operator, as in
(\ref{innerprod}) and Proposition \ref{Corth}
(here $A=\Z[1/2]$).
Similarly, the 
$\DS_w(X)$ are stable under the inclusion $\wt{W}_m \hra \wt{W}_n$,
and the products $\{\wh{\DS}_{\om}(X)\DS_{\l}(X)\}$ for $\om\in S_n$, 
$\l\in\D_{n-1}$ form an orthogonal basis for $A[X]$ as a 
$A[X]^{\wt{W}_n}$-module, where $\wh{\DS}_{\om}(X):=\sum_{w\in L(\om)}
\DS_w(X)$.

We next discuss some properties special to the type $D$ double
Schubert polynomials. Fix an integer $\ell>0$, let 
$\wt{Y}_n=(y_1,\ldots,y_{n-1},-y_n)$ and define
\[
\wh{Y}=\left\{ \begin{array}{cl}
           \wt{Y}_n
           & \mathrm{ if } \,\ n = \ell \ (\mathrm{mod}\ 2), \\
             Y_n  & \mathrm{ if } \,\ n \neq \ell \ (\mathrm{mod}\ 2).
             \end{array} \right.
\]
\begin{thm}[Maximal Grassmannian II]
\label{mgrasD2}
For each strict partition $\l\in\D_{n-1}$ of length $\ell$,
the double Schubert polynomial $\DS_{\l}(X,Y)$ is equal to 
\[
(-1)^{e(\l)+|\l'|}
\sum_{\a}\wt{P}_{\a}(X)
\sum_{\g}
(-1)^{e(\a,\g)+|\g|}
\wt{P}_{(\a\cup \g)'}(\wh{Y})
\det(e_{\g_i-\l'_j}(Y_{n-1-\l'_j})),
\]
where the first sum is over all $\a\in \D_{n-1}$ and the
second over $\g\in \D_{n-1}$ with $\g\supset \l'$,
$\ell({\g})=\ell(\l')$ and $\a\cap\g=\emptyset$.
\end{thm}
\begin{proof} If $n \neq \ell \ (\mathrm{mod}\ 2)$, then the claim
follows directly from Theorem \ref{mgrasD}. Assume that 
$n = \ell \ (\mathrm{mod}\ 2)$, fix a partition 
$\a\in\D_{n-1}$ and 
equate the coefficients of $\wt{P}_{\a}(X)$ in the sums that
occur in Theorems \ref{mgrasD} and \ref{mgrasD2}:
\begin{gather}
\label{equatel}
(-1)^{\ell}\sum_{\b}(-1)^{f(\a,\b)+|\b|}
\wt{P}_{(\a\cup \b)'}(Y_n)
\det(e_{\b_i-\l'_j}(Y_{n-1-\l'_j})) = \\
\label{equater}
\sum_{\g}
(-1)^{e(\a,\g)+|\g|}
\wt{P}_{(\a\cup \g)'}(\wt{Y}_n)
\det(e_{\g_i-\l'_j}(Y_{n-1-\l'_j})).
\end{gather}
To prove this equality, we expand each determinant 
$\det(e_{\b_i-\l'_j}(Y_{n-1-\l'_j}))$ in
(\ref{equatel}) along the last ($k$th) column, and compare with 
(\ref{equater}). The result then follows by using the identity
in the next Proposition (for varying $\m$ and $r$).
\begin{prop}
\label{keylemma}
For each partition $\m\in \D_{n-1}$ of length $r$,
we have
\[
\sum_{i=1}^r(-1)^{i-1}\wt{P}_{\m\ssm\m_i}(Y_n)\,
e_{\m_i}(Y_{n-1})=
(-1)^{r+1}\wt{P}_{\m}(Y_n)+\wt{P}_{\m}(\wt{Y}_n) .
\]
\end{prop}
\noindent
The proof of Proposition \ref{keylemma}, while elementary,
uses additional algebraic formalism, and is given in \cite[Appendix]{KT2}.
\end{proof}

\begin{cor}
For each strict partition $\l\in\D_{n-1}$ of length $\ell$,
the maximal isotropic degeneracy locus $\X_{\l}$ of 
{\em(\ref{maxlocD})} satisfies
\[
[\X_{\l}]=
(-1)^{e(\l)+|\l'|}
\sum_{\a}\wt{P}_{\a}(E^*)
\sum_{\g}
(-1)^{e(\a,\g)+|\g|}
\wt{P}_{(\a\cup \g)'}(\wh{F}^*)
\det(e_{\g_i-\l'_j}(F^{*}_{n-1-\l'_j})),
\]
where the first sum is over all $\a\in \D_{n-1}$, the
second is over $\g\in \D_{n-1}$ with $\g\supset \l'$,
$\ell({\g})=\ell(\l')$ and $\a\cap\g=\emptyset$, while
\[
\wh{F}=\left\{ \begin{array}{cl}
           \wt{F}_n
           & \mathrm{ if } \,\ n = \ell \ (\mathrm{mod}\ 2), \\
             F_n  & \mathrm{ if } \,\ n \neq \ell \ (\mathrm{mod}\ 2).
             \end{array} \right.
\]
\end{cor}

\noindent
{\bf Examples.} 1) We have
\[
\DS_{s_{\Box}}(X,Y)=
(1/2)(x_1+\ldots+x_n+y_1+\ldots+
y_{n-1}\pm y_n),
\]
with the sign of $y_n$ positive (resp.\ negative) if $n$ is even
(resp.\ odd). Also,
\[
\DS_{s_1}(X,Y)=
(1/2)(-x_1+\ldots+x_n+y_1+\ldots+
y_{n-1}\mp y_n),
\]
with the opposite sign convention for $y_n$. Note that
$\DS_{s_{\Box}}+\DS_{s_1}=\BS_{s_1}=\CS_{s_1}$ while 
$\DS_{s_i}=\BS_{s_i}=\CS_{s_i}$ for $i>1$.

\medskip
\noindent
2) Theorem \ref{mgrasD2} gives
\[
\DS_{\rho_{\ell}}(X,Y)=\sum_{\a\in \D_{\ell}}\wt{P}_{\a}(X)
\wt{P}_{\rho_{\ell}\ssm\a}(\wh{Y}).
\]

It is clear from the above examples that the type $D$ double Schubert
polynomials do not 
satisfy the same stability property as in types $B$ and $C$. However,
we see that for each $\ell>0$, the polynomial 
$\DS_{\rho_{\ell}}(X,\wh{Y})$ is stable, in the sense of \S \ref{stab}.

\section{Example: the Lagrangian Quot scheme $LQ_1(2,4)$}
\label{example}

\noindent
In this section we study the problem of extending the formula
for degeneracy loci from Corollary \ref{lagcor} to degeneracy loci 
of morphisms of vector bundles satisfying isotropicity conditions,
in analogy with the work of Kempf and Laksov \cite{KL} in type $A$.
We provide an example showing that a direct analogue of the 
Kempf-Laksov result fails in type $C$. 
This example hinges on two ingredients: Quot schemes and
degeneracy loci for a {\em morphism} from a flagged symplectic
vector bundle to a vector bundle.
Both of these, while classical for type $A$, have not
received attention in the other Lie types.
We begin with a description of degeneracy loci (we shall see that there
are two reasonable definitions, although our example will be for a
particular $\lambda\in \D_n$ for which the two definitions agree),
and later introduce the Lagrangian Quot scheme $LQ_1(2,4)$
which serves as a compactification of
the moduli space of degree $1$ maps $\bP^1\to LG(2,4)$.

\subsection{Lagrangian degeneracy loci for isotropic morphisms}
\label{lagdl}
Let $\X$ be an algebraic variety over any ground field.
Let $V$ be a symplectic vector bundle of rank $2n$ over $\X$
with complete isotropic flag of subbundles $F_i$, $1\le i\le n$, and 
let $Q$ be a vector bundle of rank $n$ over $\X$. We say that a morphism 
of vector bundles $\psi\colon V\to Q$ is {\em isotropic} if
the composite
$Q^*\to V^*\to V\to Q$ is zero, where the middle map is the 
isomorphism coming from the symplectic form on $V$. 
For such an isotropic morphism $\psi$
and for $\lambda\in \D_n$, we define the {\em Lagrangian degeneracy loci} 
$\X'_\lambda$ and $\X''_\lambda$ by
\begin{align}
\X'_\lambda &= \{\,a\in
\X\ |\ \rk(F_{n+1-\lambda_i}(a)\stackrel{\psi}\lra
Q(a))\lequ n+1-i-\lambda_i
\text{ for $1\lequ i\lequ \ell(\lambda)$}\,\}\label{ea} \\
\X''_\lambda &= \{\,a\in
\X\ |\ \rk(F^\perp_{n+1-\lambda_i}(a)\stackrel{\psi}\lra Q(a))\lequ n-i
\text{ for $1\lequ i\lequ \ell(\lambda)$}\,\}\label{eb}.
\end{align}
When $\psi$ is surjective, conditions (\ref{ea}) and (\ref{eb})
are equivalent to (\ref{Cdef}) (with $E=\Ker\,\psi$), and
$\X'_\lambda=\X''_\lambda$.
If $\psi$ is not everywhere of full rank, we only have
$\X'_\lambda\subset \X''_\lambda$ in general.
However, when $\lambda=\rho_k$ for some $k$,
the two definitions above yield the same scheme,
and we may speak without ambiguity of the Lagrangian degeneracy locus
of the morphism $\psi$. For instance, when $\l=(1)$, both 
(\ref{ea}) and (\ref{eb}) are the same as a type $A$ degeneracy
locus, so the Kempf-Laksov formula (\ref{kl1}) dictates
$[\X_1']=c_1(Q)-c_1(F_n)$.

\begin{prop}
\label{nopoly}
Fix $n=2$ and $\l=(2,1)$, and consider a general smooth variety $\X$ with 
vector bundles $V$, $Q$, $F_1$, $F_2$ and isotropic morphism
$\psi\colon V\ra Q$ as above. Then there is no polynomial in the Chern classes
of these bundles whose value equals $[\X_{\l}']$ in $CH^*(\X)$
whenever $\X_{\l}'$ has codimension $3$.
\end{prop}

We will work over $\C$ and consider the special case
where the ambient bundle $V$ as well as 
the flag of isotropic subbundles $F_{\bullet}$ are trivial. We
therefore have $V=\cO_\X\otimes W$ and $F_i=\cO_\X\otimes W_i$
for some symplectic vector space $W$ of dimension $N=2n$, with
a fixed flag of isotropic subspaces $W_1\subset\cdots\subset W_n$.
If $\psi$ is surjective (and $E=\Ker\,\psi$ as above) then for any 
$\l\in \D_n$, only the leading term of the
Lagrangian degeneracy locus formula in Corollary \ref{lagcor}
survives, hence
\begin{equation}
\label{lagtopterm}
[\X_\lambda]=\widetilde{Q}_\lambda(E^*)
\end{equation}
in $CH^*(\X)$,
provided the degeneracy locus $\X_\lambda$ of (\ref{Cdef})
has the expected dimension.

It is straightforward to give examples of the failure of (\ref{lagtopterm})
for isotropic morphisms of vector
bundles.  For instance, consider $\X=\bP^3$, 
$\l=(2,1)$, $Q=\cO(1)\oplus\cO(1)$ and the morphism
\begin{equation}
\label{abcd}
\psi \colon  \cO_{\X}\otimes W\stackrel{ \begin{pmatrix}
a&b&c&d\\
0&0&-b&a\end{pmatrix}}
\vvvvvlra Q.
\end{equation}
Here $W=\C^4$ with the standard symplectic form
\begin{equation}
\label{stdform}
\langle u,v\rangle = u_1v_3-u_3v_1+u_2v_4-u_4v_2
\end{equation}
and $[a:b:c:d]$ are the homogeneous coordinates on
$\bP^3$. To prove Proposition \ref{nopoly},
we will produce an example of an isotropic morphism of vector bundles
$V^*\to E^*$ on a nonsingular projective variety
$\X$, surjective at the generic point, 
such that the degeneracy locus $\X'_{\l}$
is smooth of the expected codimension, but whose fundamental class is 
not equal to {\em any} 
polynomial in the Chern classes of $E$.
The example was motivated by the theory of quantum cohomology,
and we explain this connection next.

\subsection{Quantum cohomology of Grassmannians}
\label{qcoh}
The theory of degeneracy loci for morphisms of vector bundles has
a direct application to the study of the quantum cohomology of 
flag manifolds in type $A$ (see e.g.\ \cite{bertram}
\cite[App.\ J]{FP} \cite{C-F}). 
In the case of the $SL_N$-Grassmannian $G$, 
Bertram \cite{bertram} used the formula 
(\ref{kl1}) of Kempf-Laksov \cite{KL} to prove a `quantum Giambelli 
formula'. The quantum Giambelli formula 
calculates the class of a Schubert variety in 
the (small) quantum cohomology ring $QH^*(G)$, with respect to 
a given presentation of this ring in terms of generators and
relations.

It is natural to copy Bertram's arguments and work towards
a quantum Giambelli formula for the type $C$ (as well as type $D$)
maximal Grassmannian varieties, contingent upon having a formula
such as (\ref{lagtopterm}) in the situation of a morphism $V^*\to E^*$, which
is generically the projection to a Lagrangian quotient bundle.
Unfortunately, this approach would dictate the {\em wrong} quantum
Giambelli formula. The simplest example is provided by the
Lagrangian Grassmannian $LG=LG(2,4)$ of isotropic $2$-planes in $\C^4$.

The (small) quantum cohomology ring $QH^*(LG)$ is generated by the special
Schubert classes $\sigma_1$ and $\sigma_2$, together with a formal variable
$q$ of degree $3$. 
The structure constants of $QH^*(LG)$ are the numbers of rational curves
on $LG$ satisfying incidence conditions (also known as Gromov-Witten
invariants).
Lines on $LG$ are parametrized by the $\bP^3$ of linear subspaces 
$\ell$ of $\C^4$:
$$\ell\subset\C^4 \quad \longleftrightarrow \quad
\{\ \Sigma\ |\ \ell\subset\Sigma\subset\ell^\perp\ \}.$$
Fixing the line in $LG$ determined by $\ell_0\subset\C^4$, as well as
a point $\Sigma_0\in LG$ in general position, i.e., satisfying
$\ell_0\not\subset\Sigma_0$, we ask how many lines on $LG$ are incident
to this line and this point.
Clearly we have conditions $\ell\subset \Sigma_0$ for the line
corresponding to $\ell$ to pass through
the point, and $\ell\subset \ell_0^\perp$ for the line to meet the given line;
thus there is a unique line on $LG$
incident to a line and a point in general position.
Consequently, the quantum product of the Schubert classes $\sigma_1$
and $\sigma_2$ receives a quantum correction term with coefficient $1$:
\begin{equation}
\label{qterm}
\sigma_1 * \sigma_2 = \sigma_{2,1} + q
\end{equation}
(here $\sigma_{2,1}$ is dual to the class of a point in $LG$ and equals
the classical cup product $\sigma_1\cup\sigma_2$).
Using the arguments of \cite{bertram}, a formula such as
(\ref{lagtopterm}) which is valid for morphisms would predict
{\em no} quantum correction term in (\ref{qterm}).
This discrepancy leads to a counterexample to (\ref{lagtopterm}) for
isotropic morphisms, which we study in detail for the remainder
of this section.

The loci $\X'_\lambda$ of (\ref{ea})
on {\em Lagrangian Quot schemes} (defined below in the situation
that we require) are the `correct' ones
for analysis of the quantum cohomology of Lagrangian Grassmannians in
analogy with \cite{bertram}. Our particular analysis will involve
$\lambda=\rho_2$, for which the loci $\X'_\lambda$, $\X''_\lambda$
coincide anyway. For more information on the quantum cohomology of
Lagrangian Grassmannians, including the full quantum Giambelli
formula for $LG(n,2n)$, see \cite{KT}.

\subsection{Quot schemes: a review}
Grothendieck's Quot schemes \cite{grothendieck}
parametrize quotients with given Hilbert polynomial
of a fixed coherent sheaf on an algebraic variety.
If $G(m,N)$ denotes the Grassmannian of $m$-dimensional subspaces of
$W
\cong \C^N$,
then a morphism $\bP^1\to G(m,N)$ is equivalent to
a map $\cO_{\bP^1}\otimes W\to Q$ to a rank $n:=N-m$ quotient bundle $Q$.
The Quot scheme $Q_d$ parametrizing
quotient {\em sheaves} of $\cO_{\bP^1}\otimes W$ with Hilbert polynomial
$nt+n+d$ is a smooth projective variety which compactifies the
moduli space of degree $d$ maps $\bP^1\to G(m,N)$.
On $Q_d\times \bP^1$ there is a universal quotient map
$\psi\colon \cO\otimes W\to T_d$.
While $T_d$ is not locally free in general,
the kernel $S_d$ of $\psi$ is locally free.
The intersection-theoretic ingredient in \cite{bertram} is the
Kempf-Laksov formula (\ref{kl1}) applied to the (nonsurjective) dual morphism
\[
\cO_{Q_d\times\bP^1}\otimes W^*\ra S^*_d.
\]

On our way to defining the Lagrangian Quot scheme,
we review the $m=1$, $N=2$ case of Quot scheme just described,
namely the Quot scheme compactification of the parameter space
$PGL_2$ of maps $\bP^1\to \bP^1$.
Let $W$ be a vector space of dimension $2$.
The Quot scheme parametrizes short exact sequences
$$0\lra S\lra W\otimes \cO_{\bP^1}\lra Q\lra 0$$
with $Q$ of rank $1$ and degree $1$.
So $S$ must be isomorphic to $\cO_{\bP^1}(-1)$, and
the Quot scheme $Q_1$ is the space of nontrivial maps of bundles
$\cO_{\bP^1}(-1)\to W\otimes \cO_{\bP^1}$,
up to multiplication by a global scalar.

Choosing a basis for $W$ and a basis $\{x,y\}$ of
$\Hom(\cO(-1),\cO)$, we have $Q_1\cong \bP^3$ with universal sheaf sequence
$$0\longrightarrow \cO(-1,-1)\stackrel{
\begin{pmatrix} ax+by\\ cx+dy\end{pmatrix}} \vlra
\cO\oplus \cO\longrightarrow Q\longrightarrow 0$$
on $\bP^3\times\bP^1$, where $[a:b:c:d]$ are the homogeneous coordinates on
$\bP^3$.

\subsection{Quot scheme of maps to $LG(2,4)$}
\label{lquot}
We now take $W=\C^4$, endowed with
the standard symplectic form (\ref{stdform}).
We introduce
the {\em Lagrangian Quot scheme} $LQ_1=LQ_1(2,4)$.
Consider the functor which parametrizes
rank $2$ degree $1$ quotients of
$\cO_{\bP^1}\otimes W$ which are (generically) Lagrangian:
the symplectic form defines an isomorphism $\cO\otimes W\to \cO\otimes W^*$,
and the Lagrangian condition on an exact sequence
\begin{equation}
\label{lagrexaseq}
0\lra S\lra \cO_{\bP^1}\otimes W\lra Q\lra 0
\end{equation}
is that the composite
$$S\lra \cO\otimes W\lra \cO\otimes W^*\lra S^*$$
is the zero map.
Since the Lagrangian condition is a closed condition,
such quotients are parametrized by a closed subscheme $LQ_1$ of the usual
Quot scheme $Q_1$.

A typical affine chart of the Lagrangian Quot scheme looks like
$$0\longrightarrow \cO\oplus \cO(-1)
\stackrel{ \begin{pmatrix}
1 & 0 \\
q & x+dy \\
r & ex+fy \\
s & gx+hy
\end{pmatrix}}\vvvlra \cO\otimes W\longrightarrow Q\longrightarrow 0$$
on $(\Spec k[d,e,f,g,h,q,r,s]/(e+qg-s, f+qh-sd))\times\bP^1$.
Such charts cover $LQ_1$, hence $LQ_1$ is smooth.

When $S\to \cO_{\bP^1}\otimes W$ is of full rank,
there corresponds a morphism from $\bP^1$ to the Lagrangian Grassmannian
$LG(2,4)$.
Since $LG(2,4)$ is a quadric hypersurface in a four-dimensional
projective space, and since lines on
$LG(2,4)$ are parametrized by $\bP^3$, we expect the parameter space
of degree $1$ maps $\bP^1\to LG(2,4)$ to be a $PGL_2$-bundle
over $\bP^3$, compactified by the Lagrangian Quot scheme.

Since it is needed below, we record the sheaf sequence corresponding
to the compactification of the space of degree 1 maps $\bP^1\to LG(2,4)$
whose image passes through a fixed point in $LG$.
The parameter space is a copy of $\bP^3\times \bP^1$,
which we endow with homogeneous coordinates $[a:b:c:d]$ and $[s:t]$
on respective factors. The sheaf sequence is the following
sequence on $\bP^3\times \bP^1\times \bP^1$:
\begin{equation}
\label{linesthroughpoint}
0\longrightarrow \cO(0,-1,0)\oplus\cO(-1,-1,-1)
\stackrel{ \begin{pmatrix}
s&0 \\
t&s(ax+by) \\
0&-t(cx+dy) \\
0&s(cx+dy)
\end{pmatrix}}\vvvvlra \cO\otimes W\longrightarrow Q\longrightarrow 0.
\end{equation}
We note the geometry: lines through a fixed point on $LG$ sweep
out a quadric cone which is a singular hyperplane section of $LG$ under its
fundamental (Pl\"ucker) embedding.
The $[s:t]$ coordinates select the line.
For each fixed $[s:t]$, the compactification of the parametrized maps
to this line is a copy of $\bP^3$, just as in the previous subsection.

\begin{prop}
\label{lquotprop}
The Lagrangian Quot scheme $LQ_1(2,4)$ compactifying the space of
degree $1$ maps from $\bP^1$ to $LG(2,4)$ is isomorphic to
the projectivization of the bundle $R\oplus R$ on $\bP^3$,
where $R$ is a rank $2$ vector bundle which fits into an exact sequence
\begin{equation}
\label{vsequence}
0\lra \cO_{\bP^3}(-1)\lra \Omega^1_{\bP^3}(1)\lra R\lra 0.
\end{equation}
\end{prop}

More naturally, in fact,
$LQ_1\simeq P(R\otimes\Hom(\cO_{\bP^1}(-1),\cO))$. Recall from Section
\ref{qcoh} that
lines on $LG(2,4)$ are parametrized by the $\bP^3$ of 1-dimensional
subspaces $\ell\subset W$.  To each such $\ell$ we associate the set of
isotropic $2$-dimensional subspaces of $W$ containing $\ell$,
or equivalently, the set of $2$-planes containing $\ell$ and
contained in $\ell^\perp$.
Thus any map $\bP^1\to LG(2,4)$ determines a one-dimensional
subspace $\ell\subset W$.

We describe what happens on the level of points, and then we
prove Proposition \ref{lquotprop} by phrasing everything in terms
of universal bundles.
In any exact sequence (\ref{lagrexaseq}),
the splitting type of $S$ must be $\cO\oplus \cO(-1)$.
As in the previous subsection, it is easier to describe the
vector bundle $S$ with injective morphism of its sheaf of
sections to $\cO\otimes W$, than it is to describe the quotient sheaf $Q$.
The line $\ell$ associated to a morphism $\bP^1\to LG(2,4)$
can be recovered from the sheaf sequence (\ref{lagrexaseq}) as
the image of the map on global sections
$$\Gamma(S)\lra W=\Gamma(\cO_{\bP^1}\otimes W).$$
This makes sense even if $S\to \cO\otimes W$ is not everywhere of full rank.
So to every point of $LQ_1$ we can associate a line in $W$.

The bundle $R$ is constructed to have fiber $\ell^\perp/\ell$
at the point of $\bP^3$ corresponding to $\ell\subset W$.
An element of $(\ell^\perp/\ell)\otimes \Hom(\cO_{\bP^1}(-1),\cO))$
determines a sheaf morphism
\begin{equation}
\label{seqa}
\cO_{\bP^1}(-1)\lra \cO_{\bP^1}\otimes(\ell^\perp/\ell).
\end{equation}
What we seek is a morphism
\begin{equation}
\label{seqb}
\cO_{\bP^1}\oplus \cO_{\bP^1}(-1)\lra \cO_{\bP^1}\otimes W
\end{equation}
with image contained in $\cO_{\bP^1}\otimes\ell^\perp$ and
whose map on global sections has image $\ell$.
Points of the projectivization of
$(\ell^\perp/\ell)\otimes \Hom(\cO_{\bP^1}(-1),\cO))$
correspond exactly to morphisms (\ref{seqa}) up to global automorphism
of $\cO(-1)$, i.e., a global scale factor.
Such a morphism (\ref{seqa}) in turn determines, uniquely, a morphism
(\ref{seqb}) which sends global sections into $\ell$,
up to global automorphism of $\cO\oplus\cO(-1)$.

\begin{proof}[Proof of Proposition \ref{lquotprop}]
Let $U\simeq\cO(-1)$ be the universal subbundle on $\bP^3$,
and let $R=U^\perp/U$, where $\perp$ refers to the standard
skew-symmetric form
$\langle\,\,,\,\,\rangle$ on $W$, and hence also on
$\cO_{\bP^3}\otimes W$.
The first claim is that $R$ fits into an exact sequence (\ref{vsequence}).
This is clear, since the (dualized, twisted) Euler sequence on $\bP^3$
$$0\lra \Omega^1(1)\lra \cO^4\lra U^*\lra 0$$
identifies $\Omega^1(1)$ with $U^\perp$.

On $\bP^3\times\bP^1$ there is the natural map of vector bundles
$$R\otimes \Hom(\cO_{\bP^1}(-1),\cO_{\bP^1})\otimes \cO_{\bP^1}(-1)\lra R.$$
Let
$P:=P(R\otimes \Hom(\cO_{\bP^1}(-1),\cO_{\bP^1}))\to\bP^3$
be the projectivization of the indicated vector bundle.
We follow the convention of \cite{F4}:
points of $P$ are points of the base together with
one-dimensional subspaces 
of the fiber. Let $S=\cO_P(-1)$ 
be the universal subbundle
of $R\otimes \Hom(\cO(-1),\cO))$ on $P$.
Then we have a morphism of bundles on $P\times \bP^1$
$$S\otimes \cO_{\bP^1}(-1)\lra R=U^\perp/U\subset (\cO\otimes W)/U.$$
Over any affine open subset of $\bP^3$, this lifts to a morphism
\begin{equation}
\label{lift}
\cO\oplus(S\otimes\cO_{\bP^1}(-1))\lra \cO\otimes W,
\end{equation}
such that, fiberwise, the map on global sections has image $U$.
To phrase this mathematically, let $\varphi$ denote the projection
$P\times\bP^1\to P$, and equally, the restriction over an affine open
of $\bP^3$.
Then $\varphi_*$ applied to (\ref{lift}) should be an isomorphism onto
$U\subset \cO\otimes W$.
This requirement determines the lift, uniquely up to automorphisms
of the fibers of $\cO\oplus(S\otimes\cO_{\bP^1}(-1))$.

Since the Quot functor is a sheaf,
these maps, defined locally,
patch to give a vector bundle $E$ and injective
map of sheaves $E\to \cO\otimes W$
on $P\times\bP^1$, such that the quotient sheaf is flat over
$P$ and on every fiber of $P\times\bP^1\to P$ has rank 1 and degree 1.
The resulting morphism $P\to LQ_1$ is an isomorphism since
it is a birational morphism of smooth varieties, which is
one-to-one on geometric points.
\end{proof}

\medskip

The natural map $\varphi^*\varphi_*E\to E$ identifies a rank 1 subbundle of $E$
with the pullback to $P\times\bP^1$ of the universal subbundle on $\bP^3$.
If we let $\pi$ denote the composite projection
$P\times \bP^1\to P\to \bP^3$, then
the cokernel of $\pi^*U\to E$,
restricted to $\pi^{-1}(A)$ for any affine open set $A\subset\bP^3$,
is identified with $S\otimes \cO_{\bP^1}(-1)$.
This identification is unique up to scaling by an
invertible regular function on $\pi^{-1}(A)$.
Hence, there is an exact sequence
$$0\lra \pi^*U\lra E\lra S\otimes \pi^*L\otimes \cO_{\bP^1}(-1)\lra 0$$
for some line bundle $L$ on $\bP^3$.
For any line $\Lambda$ in $\bP^3$,
the restriction of $E$ to
$\pi^{-1}(\Lambda)$
must be the subsheaf in the sequence (\ref{linesthroughpoint}),
up to isomorphism,
so we deduce that
$L=\cO_{\bP^3}(-1)$.
Thus we have

\begin{cor}
\label{Eseq}
Under the isomorphism $P(R\oplus R)\simeq LQ_1$ from
Proposition \ref{lquotprop}, the universal subsheaf $E$ of
$\cO\otimes W$ on $LQ_1\times\bP^1$
fits into an exact sequence
$$0\lra\cO_{\bP^3}(-1)\lra E\lra S\otimes\cO_{\bP^3}(-1)\otimes\cO_{\bP^1}(-1)
\lra 0,$$
where $S$ is the universal subbundle on $P(R\oplus R)$.
\end{cor}

\subsection{Chern class computations on $LQ_1$} \label{chernclass}
Consider, on $\X:=LQ_1\times\bP^1$, the
trivial rank $4$ bundle $V=\cO_\X\otimes W$, with standard symplectic form 
(\ref{stdform}) and
flag of isotropic subbundles $F_{\bullet}=\{F_i\}$ 
(where $F_i$ is the subbundle of
sections with nonzero entries
only among the first $i$ coordinates).
The kernel of the universal quotient map on $\X$ is the vector bundle
$E$ from Section \ref{lquot}.
Dualizing, we get a morphism
\begin{equation}
\label{esub}
\psi \colon  V^*\to E^*
\end{equation}
which is isotropic.
On $LQ_1\simeq P(R\oplus R)$ we let $h$ denote the hyperplane class on
the $\bP^3$ base of the projectivized vector bundle, and we let $z$
denote the first Chern class of $\cO_{P(R\oplus R)}(1)$
(so $z=-c_1(S)$, where $S$ is the universal subbundle on $P(R\oplus R)$).
Let $p$ be the class of a point on $\bP^1$.
Since $R$ has Chern polynomial $1+h^2t^2$, the Chern polynomial of
$R\oplus R$ is $1+2h^2t^2$.
Hence, the Chow ring (or cohomology ring) of $\X$ is given by
$$CH(\X)=\Z[h,z,p]/(h^4, z^4+2h^2z^2, p^2).$$
By Corollary \ref{Eseq}, we have
\begin{align*}
c_1(E^*) &= 2h+z+p, \\
c_2(E^*) &= h^2+hz+hp, \\
c_1(E^*)c_2(E^*) &= 2h^3+3h^2z+3h^2p+hz^2+2hzp.
\end{align*}

Recall that we have defined Lagrangian degeneracy loci
$\X'_\lambda$, $\X''_\lambda$ associated to an isotropic morphism
to a flagged vector bundle.
We consider $\lambda = (2,1)$ and the morphism (\ref{esub});
since $\lambda=\rho_2$, we have $\X'_\lambda=\X''_\lambda$.
Translating the condition (\ref{ea}) into conditions on $E$, we have
$$\X'_\lambda=\{\,a\in \X\ |\ \rk ( E(a)\ra (V/F_2)(a)) = 0\,\}.$$
If $L\subset \bP^3$ is the line which parametrizes one-dimensional
subspaces of $F_2$, then
$\X'_\lambda$ is the hypersurface in
$\pi^{-1}(L)\simeq \bP^3\times\bP^1\times\bP^1$
(where $\pi$ is the projection $\X\to \bP^3$ as in
the previous subsection) defined by the homogeneous equation
$cx+dy=0$ in the coordinates of (\ref{linesthroughpoint}).
So $\X'_\lambda$ has  
the expected codimension $3$ (and, in fact, is smooth).
Its class in the Chow ring is
\begin{equation}
\label{itsclass}
[\X'_\lambda]=h^2z+h^2p\ne c_1(E^*)c_2(E^*)=\wt{Q}_{2,1}(E^*).
\end{equation}

The above continues to be an inequality upon restriction
to $LQ_1\times \{[1:0]\}$. In fact, restricting (\ref{esub}) 
to a section of
$LQ_1\times\{[1:0]\}\simeq P(R\oplus R)\to\bP^3$ reproduces
the example (\ref{abcd}) of Section \ref{lagdl}.
Observe now that the class (\ref{itsclass}) is not
equal to any polynomial in the Chern classes of $E$.
In particular, this proves Proposition \ref{nopoly}.

\end{document}